\documentclass{amsart}
\usepackage{amsthm}
\newtheorem{lemma}{Lemma}

\newtheorem{theorem}{Theorem}
\newtheorem{corollary}{Corollary}
\newtheorem{definition}{Definition}

\newtheorem{remark}{Remark}

\usepackage{changebar}

\newtheorem{fact}{Fact}

\newcommand{\nmid}{\not \hspace{0.25em} \mid}
\newcommand{\nth}[1]{$#1 {\rm - th }$}

\newcommand{\Z}{\mathbb{Z}}
\newcommand{\ZM}[1]{\Z /( #1 \cdot \Z)}
\newcommand{\ZMs}[1]{(\Z / #1 \cdot \Z)^*}

\newcommand{\ord}{{\rm ord}}
\newcommand{\rem}{{\rm  rem }}

\newcommand{\Tr}{\mbox{\bf Tr}}

\newcommand{\rg}[1]{\mbox{\bf #1}}
\newcommand{\eu}[1]{\mathfrak{#1}}
\newcommand{\id}[1]{\mathcal{#1}}

\newcommand{\GCD}{\mbox{ GCD }}
\newcommand{\Aut}{\mbox{ Aut }}
\newcommand{\End}{\mbox{ End}}
\newcommand{\GL}{\mbox{ GL}}
\newcommand{\rf}[1]{(\ref{#1})}

\newcommand{\Norm}{\mbox{\bf N}}
 
\newcommand{\lchooses}[2]{\left( \frac{#1}{#2 } \right)}

\newcommand{\F}{\mathbb{F}}
\newcommand{\K}{\mathbb{K}}

\newcommand{\Q}{\mathbb{Q}}

\newcommand{\N}{\mathbb{N}}
\newcommand{\KH}{\mathbb{H}}

\def\NI{\noindent}

\newcommand {\legendre}[2]{\genfrac {(}{)}{1pt}{}{#1}{#2}}

\begin{document}

\title[Dual Elliptic Primes]{Dual Elliptic Primes and Applications to Cyclotomy Primality
Proving} 
\author{Preda Mih\u{a}ilescu}
\address{Mathematisches Institut der Universit\"at
G\"ottingen}
\email{preda@uni-amth.gwdg.de} 

\thanks{The research was completed while the author is holding a
research chair sponsored by the Volkswagen Stiftung} 
\date{Version 2.0 \today}
\bigskip

{\obeylines \em 
\vspace*{0.5cm} 
\noindent\hspace*{2.0cm}Quand j'ai couru chanter ma p'tit' chanson pour Marinette \newline
\vspace*{-0.4cm}
\noindent\hspace*{2.0cm} La belle, la tr\^{a}itresse \'{e}tait all\'{e}e \`{a} l'op\'{e}ra.\newline
\vspace*{-0.4cm}
\noindent\hspace*{2.0cm} Avec ma p'tit' chanson, j'avais l'air d'un con, ma m\`{e}re.\newline
\vspace*{-0.4cm}
\noindent\hspace*{2.0cm} Avec ma p'tit' chanson, j'avais l'air d'un con\footnote{Georges Brassens: \textit{Marinette}}. \newline

\hspace*{5.8cm}To Hendrik W. Lenstra, Jr.
\vspace*{1.0cm}
}

\begin{abstract}
Two rational primes $p, q$ are called {\em dual elliptic} if there is
an elliptic curve $E \mod p$ with $q$ points. They were introduced as
an interesting means for combining the strengths of the elliptic curve
and cyclotomy primality proving algorithms. By extending to elliptic
curves some notions of galois theory of rings used in the cyclotomy
primality tests, one obtains a new algorithm which has heuristic cubic
run time and generates certificates that can be verified in quadratic time. 

After the break through of Agrawal, Kayal and Saxena has settled the
complexity theoretical problem of primality testing, some interest
remains for the practical aspect of state of the art implementable
proving algorithms.
\end{abstract}
\maketitle
\section{Introduction}
Primality testing is a discipline in which {\em constructions of objects } in
fields of positive characteristic $p$ are mimicked in algebras over rings
$\ZM{n}$ for integers $n$ which one believes to be prime, and of whose
primality one wishes to have a proof. The constructions should then allow an
efficient computation and be based on operations which have the property of
either yielding results over $\ZM{n}$ or else display a factor of $n$ or at
least a proof of its compositeness.

In the simplest cases, the constructions restrict to simple
verifications. Fermat's ``small Theorem'' stating that $a^{p-1} \equiv 1 \mod
p$ for rational primes $p$ and bases $a$ not divisible by $p$, is the first
ingredient used for fast verification of primality of integers $n$. In the
simplest version of the idea, the \textit{Fermat pseudoprime test, to base
$a$} checks $a^{n-1} \equiv 1 \bmod n$ and returns ``composite'', if the
congruence is not verified. If it is verified, only probabilistic statements
can be made about primality of $ n$.

Stronger statements are obtained when one has sufficient information about the
factorization of $n-1$. For instance, if there is a prime $q | (n-1)$ and $q >
\sqrt{n}$, while $\left(a^{(n-1)/q}-1,n\right) = 1$ and $a^{n-1} \equiv 1
\bmod n$, then one easily proves that $n$ is prime. This test constructs a
{\em primitive \nth{q} root } of unity modulo $n$, in the sense that
$\Phi_q(\alpha) = 0 \bmod n$ with $\alpha = a^{(n-1)/q} \rem \ \ n$ and
$\Phi_q(X)$ the \nth{q} cyclotomic polynomial. Tests of this type are known
under the name of Lucas - Lehmer tests. They share the feature, that one
proves that a certain number $a \in \ZMs{n}$ is a primitive \nth{q} root of
unity for some $q > \sqrt{n}$ - so it generates a cyclic subgroup of $\ZMs{n}$
which is, by its size, incompatible with the hypothesis that $n$ is be
composite.

The idea was generalized, freeing it of the requirement for a priori knowledge
of large factors of $n-1$. This is made possible by working in larger
extensions of $\ZM{n}$ and using more involved properties of rings in
cyclotomic fields and the related Gauss and Jacobi sums. The resulting
algorithms are currently denoted by the generic name \textit{Cyclotomy
Primality Proving} (CPP). They originate in the work of Adleman, Pomerance and
Rumeley \cite{APR} and were improved by Lenstra et. al. \cite{Le1},
\cite{Le2}, \cite{Le3}, \cite{CoLe}, \cite{BH}, \cite{Mi1}. Their main idea is
to building a frame -- a Galois algebra over $\ZM{n}$ -- in which a factor
$\Psi(X) | \Phi_s(X) \bmod n$ can be constructed for some large $s$ and such
that, if $n$ is prime, the factor is irreducible. The definitions of the
Galois algebras in which the test take place have undergone some variations
\cite{BH, Mi1, MV, LP} since their introduction in \cite{Le2}.

The name CPP covers an unconditionally deterministic variant and one which is
deterministic under assumption of the ERH, as well as a Jacobi sum and a Lucas
- Lehmer variant; all the variants may well be combined together. The CPP test
provides a proof of the fact that the \nth{s} cyclotomic polynomial $\Phi_s(X)
\in \Z[X]$ -- for some special, large and highly composite integers $s$ --
factors modulo $n$ the way it should, if $n$ were prime. If this is the case,
primality of $n$ follows, \textit{or } the existence of some prime factor
\begin{eqnarray}
\label{findiv}
  r \in \{ n^i \ \rem \ s : i = 1, 2, \ldots, t = \ord_s(n) \} .
\end{eqnarray}

The algorithms of CPP are \textit{de facto} fast, competitive
primality proving algorithms, but they have the complexity
theoretical intolerable feature of a provable {\em superpolynomial}
run - time 
\begin{eqnarray}
\label{spol}
O\left(\log(n)^{\log \log \log(n)}\right),
\end{eqnarray}
which is in fact the expected size of $t$ in \rf{findiv}.

The use of elliptic curves was first proposed for primality proving by
Goldwasser and Kilian \cite{GoKi} in an algorithm which was proved to be
random polynomial up to a possible, exponentially thin, exceptional set. The
algorithm was made computationally practical by Atkin \cite{AtMo} who
suggested a method of determining the expected number of points on an elliptic
curve, by using complex multiplication. It now runs under the generic name
ECPP (Elliptic Curve Primality Proving) and was first implemented in 1989 and
continuously improved since then, by F. Morain \cite{ECPP}.

The algorithms we present in this paper build up upon the idea of Atkin on the
one hand, on extending the use of Galois rings to the context of elliptic
curve primality proving and, finally, on a novel concept of \textit{dual
elliptic primes}. These are loose relatives of \textit{twin primes} in
imaginary quadratic extensions and allow to combine the worlds of CPP and ECPP
in a new algorithm that we call CIDE. The fundamental gain of CIDE consists in
eliminating the alternative \rf{findiv} in CPP, thus yielding a random
polynomial algorithm, which is practically an improvement of both CPP and
ECPP. We note that the computation Jacobi sums, which was an other
superpolynomial step in CPP, can be solved in random polynomial time thanks to
the novel algorithm of Ajtai et. al. \cite{AjKuSi}; in practice, the
computation of Jacobi sums can be solved in very short time using their
arithmetic properties and a PARI program for finding generators of principal
ideals. Herewith CIDE is faster by a factor of $\log(n)$ then either version
of ECPP; i.e. the \cite{GoKi}, which is slower but has a proof of random
polynomial run time for \textit{almost all} inputs, or FastECPP \cite{Mo2},
which runs de facto in time $O\left((\log(n)^{4 + \varepsilon}\right)$, but
the run time proof uses some heuristics. Unsurprisingly, the same kind of
proofs can be provided for the two versions of CIDE: this is due to the fact
that the first step of finding a pair of dual elliptic pseudoprimes requires
running one round of some version of ECPP.

The structure of the paper is the following. In the next section we
give some general definitions and facts related to elliptic curves
over finite fields, complex multiplication and ECPP. In the third
section we develop a theory of elliptic extensions of galois rings,
which is a natural analog of cyclotomic extensions used in CPP
\cite{Mi4}. Section four brings the definition of dual elliptic primes
and their pseudoprime counterparts and the basic properties of
pseudoprimes which are going to be exploited algorithmically in the
subsequent section. Finally, section six gives run time analysis and
implementation data and in section seven we draw some brief
conclusions.

\section{Elliptic curves and related pseudoprimes}
If $\K$ is some field, the equation $Y^2 \equiv X^3 + A X + B$,
with $A, B; X, Y \in \K$ and the discriminant $\Delta = 4 A^3 + 27 B^2 \neq 0$,
 defines an elliptic curve over $\K$. We denote it by
\begin{eqnarray}
\label{edef}
 \id{E}_{\K}(A, B) = \{ \ (X, Y) \in \K^2 : Y^2 = X^3 + A X + B \ \},
\end{eqnarray}
or simply $\id{E}$ when there is no ambiguity. The elements $P = (X, Y) \in
\id{E}$ are {\em points} and the curve is endowed with an addition law, $R = P
\oplus Q$ defined by
\begin{eqnarray}
\label{ede}
\lambda & = & \frac{Q_y - P_y}{Q_x - P_x}, \quad \hbox{ for } \quad P
\neq Q, \nonumber \\ 
\lambda & = & \frac{3 P_x^2 + a}{2 P_y}, \quad \hbox{ for }
\quad P = Q, \\
R_x & = & \lambda^2 - (P_x + Q_x), \quad R_y = \lambda R_x
+ (P_y - \lambda P_x). \nonumber
\end{eqnarray}
We let 
\[ \mu(P, Q) = \begin{cases} Q_x - P_x & \hbox{if} P \neq Q, \\ 2 P_y & \hbox{otehrwise} \end{cases}.\]
The neutral element is the \textit{point at infinity} $\eu{O}$ and $P \oplus Q
= \eu{O}$ iff $\mu(P, Q) = 0$; the inverse of $P = (X, Y)$ is $-P = (X, -Y)$.
This makes $\id{E}$ into an abelian group - see also \cite{W}, \S 2.2.  The
$k$ - fold addition of a point with itself is written $[ k ] P$ and can be
expressed by explicite polynomials over $\K$:
\begin{eqnarray}
\label{mult} 
[ k ] P = \left(\frac{\phi_n(P_x)}{\psi_n^2(P_x)}, P_y
\frac{\omega_n(P_x)}{\psi^3_n(P_x)} \right), \quad \hbox{with } \quad \phi_n,
\psi_n, \omega_n \in \Z[ A, B ],
\end{eqnarray}
see \cite{W}, Theorem 3.6, where the $Y$ coordinates are given by some
bivariate polynomials. These can be reduced to mono-variate ones as above.

The $k$ - torsion of $\id{E}_{\K}(A, B)$ is the set
\[   \id{E}_{\K}(A, B)[ k ] = \left\{ P \in \id{E}_{\overline{\K}}(A, B) :  [ k ] P = \eu{O} \right\}.\]
Note that the torsion if defined over the algebraic closure; if the
characteristic is $0$ or coprime to $k$, then $ \id{E}_{\K}(A, B)[ k ] \cong
\ZM{k} \oplus \ZM{k}$, e.g. \cite{W}, Chapter 3. Furthermore, the torsion is
related to the zeroes of $\psi_k(X)$ by
\begin{eqnarray}
\label{psi}
 \id{E}_{\K}(A, B)[ k ] = \left\{P \in \id{E}_{\overline \K}(A, B) \ : \
 \psi_k(P_x) = 0. \right\}
\end{eqnarray}

In algorithmic applications, the field $\K$ is a finite field. Here it is
mostly a prime field $\F_p$, with $p$ a rational prime and we write
$\id{E}_{\F_p} = \id{E}_p$.  In this case, the size of the group is bounded by
the Hasse interval
\[ m = \left| \id{E}_p \right| \in \left((\sqrt{p}-1)^2, \sqrt{p}+1)^2 
\right) .\]

It is useful to consider the addition law of elliptic curves also over rings
$\ZM{n}$, with $n$ a rational integer, which needs not be a prime. In such
cases the addition law is not everywhere defined, but it turns out that
exactly the points $P, Q$ for which $P \oplus Q$ is not defined are of great
algorithmic use. The application of this generalization are found in factoring
and primality testing. Since the conditions which are given in fields by $T
\neq 0$ -- e.g. for $T = \mu(P., Q)$ or $T = \Delta$ -- are replaced by GCD
computations and the requirement that $T \in \ZMs{n}$, whenever such a
condition is not met, a possible non trivial factor of $n$ is found. Thus the
fact that addition is not defined in such a case turns out to be an advantage
rather then a nuissance, since finding non trivial factors achieves the goal
of the algorithm.

Formally, for a given $n \in \N_{>1}$ one lets
\begin{eqnarray}
\label{endef}
 \id{E}_{n}(A, B) & = & \{ \ (X, Y) \in (\ZM{n})^2 : Y^2 = f(X) \ \}, \quad \hbox{ with} \\
  f(X) & = & X^3 + A X + B \nonumber
\end{eqnarray}
where $A, B \in \ZM{n}$ are such that $4 A^3 + 27 B^2 \in \ZMs{n}$. Addition
of two points is defined by \rf{ede} \textit{whenever} $\mu(P, Q) \in
\ZMs{n}$. Certainly, the pair $(\id{E}_n, \oplus)$ does not define a curve in
the sense of algebraic geometry and is not even a group. We may however and
shall refer to the set of points $\id{E}_n(A, B)$ as the elliptic curve with
parameters $A, B$ over $\ZM{n}$ and use the partial addition on this curve.

In primality testing we have the usual ambiguity consisting in the fact that
the curves $\id{E}_n$ which we use are defined in the sense of \rf{endef}; if
a test for $n$ completes successfully, they turn out to be proper curves in
the sense of algebraic geometry, defined over the field $\F_n$. Otherwise, non
trivial factors of $n$ or other contradictions to the hypothesis that $n$ is a
prime may be encountered in the process of a test.

Due to \rf{mult}, the $k$ - fold addition can be uniquely defined for any $P
\in \id{E}_n(A, B)$ such that $\psi_k(P_x) \in \ZMs{n}$; it does not depend on
particular addition chains for $k$. Note that since $A, B \in \ZM{n}$ and
$\psi_k \in \Z[ A, B ]$, the division polynomial $\psi_k(X) \in \ZM{n}[ X ]$.
Let the $k$ - torsion in this case be
\[   \id{E}_{n}(A, B)[ k ] = \left\{ P \in \id{E}_{n}(A, B) :  (\psi_k(P_x),
n) \neq 1 \right\}.\] We say that a torsion point $P \in \id{E}_{n}(A, B)$ is
proper, if $(\psi_k(P_x), n) = n$; for an improper $k$ - torsion point, an
algorithm using $k$ - multiplication on $\id{E}_{n}(A, B)$ would end by
featuring a non trivial divisor of $n$.

Note that unlike the field case, we have only defined torsion points of
$\id{E}_n(A, B)$ which lay in $(\ZM{n})^2$. For the general case, we need a
substitute for the algebraic closure of a field. For this we define the
following formal algebras:
\begin{definition}
\label{toralg}
Let $\rho_k(X) | \psi_k(X)$ be a polynomial such that $(\rho_k(X), \psi_i(X))
= 1$ for $i < k$. We define a \textbf{k-torsion algebra} $\rg{R}$ and the
\textit{two points} $k$-torsion algebra $\rg{R'}$ by:
\begin{eqnarray}
\label{toral}
\rg{R} & = & \ZM{n}[ X ]/(\rho_k(X)) \quad \hbox{ and } \Theta = X \bmod
\rho_k(X) \in \rg{R} , \\ \rg{R'} & = & \rg{R}[Y]/\left(Y^2 -
f(\Theta)\right), \quad \Omega = Y \bmod \left(Y^2 - f(\Theta)\right) \in
\rg{R'}. \nonumber
\end{eqnarray}
\end{definition}
In an two points torsion algebra $\rg{R}'$, the pair $P = (\Theta, \Omega) \in
\rg{R'}^2$ verifies by construction the equation of $\id{E}_n(A, B): Y^2 =
f(X)$.

We claim that the iterated addition $[ i ] P$ is defined for $P$ and each $i <
 k$. Indeed, if this were not the case for some $i < k$, there is a prime $p |
 n$ and a maximal ideal $\eu{P} \subset \rg{R}'$ containing $p$, such that $[
 i ] P \bmod \eu{P} = \eu{O}_p$, the point at infinity of the curve
 $\id{E}_{\overline{\F_p}}(A \bmod p, B \bmod p)$. This contradicts the
 premise $(\rho_k(X), \psi_i(X)) = 1$, thus confirming the claim. It follows
 that the points $[ i ] P \in \rg{R'}^2$ are $k$ - torsion points in the
 two points algebra \footnote{We are not interested here in the problem of
 constructing algebras which contain, like in the field case, two linear
 independent torsion points.}.
 
There is a unique monic polynomial $g_i(X) \in \ZM{n}[ X ]$ of degree $<
 \deg(\psi_k(X))$, such that $\psi_i^2(X) \cdot g_i(X) \equiv \phi_i(X) \bmod
 \psi_k(X)$. Then $g_i(\Theta) = ([ i ] P)_x$, by \rf{mult}, since
 $\psi_k(\Theta) = 0$. We have thus:
\begin{eqnarray}
\label{mulpol} g_i(\Theta) = ([ i ] P)_x, \quad \hbox{with} \quad P = 
(\Theta, \Omega) \in \rg{R'}^2.
\end{eqnarray}

A {\em size} $s\left(\id{E}_n\right)$ will be the result of some algorithm for
computing the number of elements of an elliptic curve in the case when $n$ is
prime. The size may depend upon the algorithm with which it is computed. Two
approaches are known: the variants of Schoof's algorithm \cite{Sch} and the
complex multiplication approach of Atkin \cite{AtMo}.

We can herewith extend some notions of pseudoprimality to elliptic curves:
\begin{definition}
Let $n$ be an integer and $\id{E}_n(A, B)$ a curve with size $m$. We say that
$n$ is \textit{elliptic Fermat pseudoprime} with respect to this curve, if
there is a point $P \in \id{E}_n(A, B) \in \id{E}_n(A, B)[ m ]$.

Furthermore, if $q | m$ is an integer, we say that $n$ passes an
\textit{elliptic Lucas - Lehmer} test of order $q$ (with respect to
$\id{E}_n(A, B)$), if there is a point $P \in \id{E}_n(A, B)[ q ]$.
\end{definition}

The test of Goldwasser and Kilian \cite{GoKi}, which is the precursor of ECPP,
can herewith be stated as follows: given $n$, find a curve $\id{E}_n(A, B)$
with a size $m$ divisible by a probable prime $q > (p^{1/4}+1)^2$ and show
that $n$ passes a Lucas - Lehmer test for $q$. If $q$ is an actual prime, then
the test implies that $n$ is also one. So one iterates the procedure for $q$,
obtaining a descending chain which reaches probable primes of polynomial size
in $O(\log(n))$ steps. In \cite{GoKi} sizes are estimated using the algorithm
of Schoof.  Even in the much faster version of these days \cite{BMSS}, this
would still yield an impractical algorithm. It does have the advantage of a
provable run time analysis.

If the field $\K = \F_q$ is a finite field of characteristic $p$, then the
Frobenius map $\Phi_q: X \mapsto X^q$ is an endomorphism of $\id{E}_{\overline
\F_q}(A, B)$ and verifies a quadratic equation:
\begin{eqnarray}
\label{Frobeq}
   \Phi_q^2 - t \Phi_q + q = 0
\end{eqnarray}
in $\End\left(\id{E}_{\overline \F_q}(A, B)\right)$, as shown for instance in
\cite{Si}, p. 135. The number $t$ is related to the size of the group $\id{E}$
over $\F_q$ by $\left| \id{E}_q \right| = q + 1 - t$. In particular, if $q = p
= \pi \cdot \overline \pi$, for $\pi \id{O} \subset \K$, the ``CM field'' of
$\id{E}$ (see below), then $t = \Tr(\pi)$, \cite{Cox} Chapter 14, in
particular Theorem 14.6 .

The Frobenius acts as a linear map on $\id{E}_{q}(A, B)[k]$. If $k = \ell$ is
a prime, $\id{E}[ \ell ]$ is a vector space and there is a matrix
$M_{\ell}(\Phi_q) \in \GL_2(\F_{\ell})$ associated to the Frobenius modulo
$\ell$.  The reduced equation \rf{Frobeq} modulo $\ell$ is also the
characteristic polynomial of $M_{\ell}(\Phi_q)$.

If $\delta = t^2 - 4 q$ is a quadratic residue over
$\F_{\ell}: \lchooses{\delta}{\ell} = 1$, then the equation
\rf{Frobeq} has two distinct roots $\bmod \ell$, which are the {\em
eigenvalues} $\lambda_{1, 2} \in \F_{\ell}^{\times}$ of the
Frobenius. Accordingly, there are points $P_{1, 2} \in \id{E}_q(A, B)[\ell]$
such that
\[ \Phi_q(P_i) = [ \lambda_i ] P_i, \quad i = 1, 2. \]

In the context of algorithms for \textit{counting points on elliptic curves }
\cite{Sch}, the primes with $\lchooses{\delta}{\ell} = 1$ are often referred
as {\em Elkies primes}, while all other primes are {\em Atkin} primes. In this
case, to each eigenvalue there corresponds an \textit{eigenpolynomial }
defined by
\begin{eqnarray}
\label{eigpol}
F_i(X) = \prod_{k=1}^{(\ell-1)/2} \left(X - \left( [ k ] P_i\right)_x \right)
\in \F_q[ X ], \quad i = 1, 2.
\end{eqnarray}
Here $\left( [ k ] P_i\right)_x$ is the $x$ - coordinate of the point $[ k ]
P_i$. Various algorithms have been developed for the fast computation of the
eigenpolynomials, without prior knowledge of the eigenpoints or eigenvalues;
see for instance \cite{BMSS} for a recent survey.

\subsection{Complex Multiplication and Atkin's approach to ECPP}
We recall some facts about complex multiplication and refer to \cite{Cox},
Chapter 14 and \cite{Si}, Chapter V, for more in depth treatment.
\begin{fact}
Let $p$ be a prime and $\id{E}_p(A, B)$ be an ordinary elliptic
curve\footnote{A curve is ordinary if it is regular and not supersingular,
\cite{W}, p. 75}. Then there is a quadratic imaginary field $\K =
\Q(\sqrt{-d})$ and an order $\id{O} \subset \K$ such that:
\begin{itemize}
\item[1.] The endomorphism ring of $\id{E}_p(A, B)$ is isomorphic to $\id{O}$.
\item[2.] There is a $\pi \in \id{O}$ such that $p = \pi \cdot \overline \pi$
  and the number of points 
\begin{eqnarray}
\label{sizAtk}
 \left| \id{E}_p(A, B) \right| = \Norm (\pi \pm 1),
\end{eqnarray}
the sign being defined only up to twists.
\item[3.] If $H_{\id{O}}(X) \in \Z[ X ]$ is a polynomial which generates the
  ring class field $\KH$ of $\id{O}$, i.e. $\KH = \K[ X ]/ \left(
  H_{\id{O}}(X) \right)$, then $H_{\id{O}}(X)$ splits completely modulo $p$.
\item[4.] There is a zero $j_0 \in \KH$ of the polynomial $H_{\id{O}}(X)$ and
an elliptic curve $\id{E}_{\KH}(a, b)$ defined over $\KH$ such that:
\begin{itemize}
\item[a)] The $j$ -invariant of $\id{E}_{\KH}(a, b)$ is $j_0$, of $r(j_0)$
with $r(X) \in \Q(X)$.
\item[b)] Its endomorphism ring is isomorphic to $\id{O}$ and
\item[c)] The curve has good reduction at a prime $\wp \subset \id{O}(\KH)$
above $(p)$.
\item[d)] The reduction is $\id{E}_p(A, B)$ and it is a direct consequence of
  CM, that $\id{E}_p(A, B)$ is ordinary.
\end{itemize}
Under these circumstances, the curve $\id{E}_{\KH}(a, b)$ is unique and is
called the \textit{Deuring lift} of $\id{E}_p(A, B)$.\end{itemize}

In $\id{O}$, the prime $p$ splits in principal ideals in $\id{O}$ if and only
if 3. holds - see e.g. \cite{Cox} Theorem 9.2 for the case when $\id{O}$ is
the maximal order. In particular:
\begin{eqnarray}
\label{quform}
 p = \pi \cdot \overline \pi \quad \hbox{ with } \quad \pi \in \id{O} \quad
\Leftrightarrow \quad \exists \ x \in \F_p: \ H_{\id{O}}(x) = 0.
\end{eqnarray}
\end{fact}

Thus the endomorphism ring associates an order in an imaginary quadratic field
to an ordinary elliptic curve over a finite field - the association being
actually an isomorphism of rings. Non-isomorphic curves can be associated to
one and the same order. This fact allows to construct curves over a finite
field $\F_p$ which have a known endomorphism ring and thus the size may be
derived directly from \rf{sizAtk}. The algorithm involves the construction of
polynomials $H_{\id{O}}(X)$ for various orders of small discriminant until one
is found which splits completely modulo $p$. The methods for computing
$H_{\id{O}}(X)$ have been subject of investigation for over a decade; see
\cite{Mo3} for an in depth treatment and \cite{BMSS} for current
improvements. The advantage of this approach, is that curves with known size
can be computed faster then by using the best versions of Schoof's algorithm
for computing the size of a given curve. Thus although this approach is not
used for finding the size of a given curve, it is sufficient for some
application where it suffices to know \textit{some} curve together with its
number of points.

The idea of Atkin was to produce similar associations for curves $\id{E}_n(A,
B)$, with $n$ not necessarily prime, and to estimate their size using the
equation in \rf{sizAtk}. In order to produce such an association, one uses
algorithms for finite fields. The construction may thus stop with a
contradiction to the hypothesis that $n$ is prime. Otherwise it is expected to
produce an order $\id{O} \subset \Q[\sqrt{-d}]$ in which $n$ factors in
principal ideals $n = \nu \cdot \overline \nu: \ \nu \in \id{O}$ and such that
$H_{\id{O}}(X)$ has a linear factor in $\ZM{n}$. Furthermore, it produces a
curve $\id{E}_n(A, B)$ with \textit{Atkin size} $m = \Norm(\nu \pm 1)$ as
suggested by \rf{sizAtk}. Several discriminants $d$ are tried, until it is
found by trial factorization that $m$ is divisible by a large pseudoprime
$q$. Finally, a point $P \in \id{E}_n(A, B)$ is sought, such that $\psi_q(P_x)
\not \in \ZMs{n}$. If $P$ is not a proper \nth{q} torsion point, a non trivial
factor of $n$ is found and the algorithm terminates. Otherwise, if $q$ is in
fact prime, then so must $n$ be, by the Lucas - Lehmer argument. This leads to
an iterative primality proof, like in the case of Goldwasser and Kilian, but
with a faster estimation of the size. However, since the discriminants $d$
must have polynomial size, the curves taken into consideration are not
random. Unlike the case of \cite{GoKi}, the fact that one can find in
polynomial time a discriminant such that the above conditions hold is
supported by heuristic arguments. Such arguments are given in \cite{GaK}.

We introduce the following notion of pseudoprimes, related to the above
algorithm:
\begin{definition}
\label{Atkassoc}
Let $n$ be an integer and $\id{E}_n(A, B)$ be an elliptic curve (with partial
addition), $\K = \Q(\sqrt{-d})$ be a quadratic imaginary field and $\id{O}
\subset \K$ some order. We say that $\left(\id{E}_n(A, B), \id{O}\right)$ are
{\em associated} if the following conditions are fulfilled:
\begin{itemize}
\item[1.] The integer $n$ is square free, there is a $\nu \in \id{O}$ such
that $n = \nu \cdot \overline \nu$, and
\begin{eqnarray}
\label{coprim}
\left(n, \nu + \overline \nu\right) = 1.
\end{eqnarray}
\item[2.] There is a polynomial $H_{\id{O}}(X) \in
\Z[X]$, which generates the ring
class field $\KH$ of $\id{O}$, i.e. $\KH = \K[ X ]/ \left(
H_{\id{O}}(X) \right)$ and which has a zero $\jmath_0 \in
\ZMs{n}$. Furthermore, the $j$ - invariant of $\id{E}_n(A, B)$ is a rational
function in $\jmath_0$. 
\end{itemize}
\end{definition}
\begin{remark}
We refer the reader to \cite{AtMo, EnMo02, EnMo03} for details on techniques for choosing
the polynomial $H_{\id{O}}$. It should be mentioned that the {\em modular
equation} is a theoretical alternative for the polynomial $H_{\id{O}}(X)$, and
it has the $j$ - invariants as zeroes; however, from a computational point of
view, the modular equation is impractical, having very large coefficients, so
one constructs alternative polynomials which generate the same field.
\end{remark}

Based on the associations of curves and orders, one defines Atkin pseudoprimes
as follows:
\begin{definition}
\label{Atkps}
We say that $n$ is \textbf{Atkin pseudoprime}, if
\begin{itemize}
\item There is a curve $\id{E}_n(A, B)$ associated to an order $\id{O} \subset
\K = \Q[\sqrt{-d}]$ according to the above definition.
\item The Atkin size of $\id{E}_n(A, B)$ is $m = \Norm(\nu \pm 1)$ and is divisible by a strong pseudoprime $q > \left(n^{1/4} + 1\right)^2$.
\item There is a proper $q$ - torsion point $P \in \id{E}_n(A,
B)$\footnote{Since $q | m$, a $q$ torsion point should be found in the curve
over $\ZM{n}$, if $n$ is prime, so the condition is consistent.}
\end{itemize}
The pseudoprime $n$ is thus given by the values 
\[ \left( n; (\id{E}_n(A, B), \id{O}); P, q \right) . \] 
\end{definition}

In all versions of the ECPP test, one seeks a random curve whose size is
divisible by some large pseudoprime $q$. When the parameters $A, B \in \ZM{n}$
are chosen uniformly random. In this case, if $n$ is a prime, it is known that
the sizes of the curves are close to uniform distributed in the Hasse interval
\cite{CP}, Theorem 7.3.2. This fact is useful for the run time analysis of the
Goldwasser - Kilian test.

Atkin's test builds descending sequences of Atkin pseudoprimes $n, q, \ldots$,
until pseudoprimes of polynomial size are reached.  The discriminant $-d$ of
the field $\K$ must be polynomial in size, which is an important restriction
for the choice of $\id{O}$. For prime $n$, the density of the curves with CM
in fields with polynomial discriminant is exponentially small. Thus Theorem
7.3.2 does not hold and there is thus no {\em proof} for the fact that ECPP
terminates in polynomial time even on {\em almost } all inputs.

We note the following consequence of condition 2.:
\begin{lemma}
\label{l1}
Suppose that $n > 2$ is an integer for which there exists an association
$(\id{E}_n(A, B), \id{O})$ according to Definition \ref{Atkassoc} and let $p |
n$ be a rational prime. Then $\id{E}_p(\overline A,\overline B)$ with
$\overline A = A \ \rem \ p, \overline B = B \ \rem \ p$ is an elliptic curve
over the field $\F_p$ with CM in $\id{O}$ and $p$ splits in principal ideals
in this order, say $p = \pi \cdot \overline \pi$.
\end{lemma}
\begin{proof}
The curve $\id{E}_p(\overline A, \overline B)$ is defined by reduction modulo
$p$. The polynomial $H_{\id{O}}(X)$ has a root $\jmath_0 \in \ZM{n}$ and thus
$\overline \jmath_0 = \jmath_0 \bmod p$ is a root thereof in $\F_p$. The $j$
invariant will then be a rational function of this value. Then \rf{quform}
implies that $p = \pi \cdot \overline{\pi}$.
\end{proof}

\section{Gauss sums and CPP}
The Jacobi sum test \cite{APR, Le1}, which is the initial version of CPP is
based on the use of Gauss and Jacobi sums. Over some field $\K$, these are
classical character sums, see e.g. \cite{IR}, Chapter 8. In primality testing
however, the images of the characters are taken over some ring $\ZM{n}$ which
need not be a field. We need thus a dedicated context of \textit{cyclotomic
extensions of rings} for the definition of these sums.

Since their definition by Lenstra~\cite{Le2}, cyclotomic extensions have
undergone various modifications \cite{BH, Mi1, Mi2, Le4} until the recent
``pseudo-fields''~\cite{Le4, LP}.  We shall follow use here definitions given
in ~\cite{AM, Mi2}. Proofs of the facts we shall need are in \cite{Mi1, Mi2}.

Let $n \in \N$ be an integer and consider rings of characteristic $n$, more
precisely finite Abelian ring extensions $\rg{R} \supset \ZM{n}$.  Galois
extensions~\cite{Mi2} are simple algebraic extensions of the form $\rg{R} =
\ZM{n}[ T ]/(f(T))$ endowed with automorphisms which fix $\ZM{n}$.  We are
interested in the \textit{simple Frobenius extensions} defined by:
\begin{definition}
\label{defgalext}
Let $\rg{R}$ be a finite commutative ring of characteristic $n$ and $\Psi(X)
\in \rg{R}[X]$ a monic polynomial. We say that the ring extension $\rg{R} =
\ZM{n}[X]/(\Psi(X))$ is simple Frobenius if:
\begin{itemize}
\item[F1.] There is a $t > 0$ such that
\[  \Psi(X) = \prod_{i=1}^t \ \left(X - \zeta^{n^i} \right), \quad 
\hbox{ where } \quad \zeta = X + (\Psi(X)) \in \rg{R}. \]
\item[F2.]  Let $x_i = \zeta^{n^i} \in \rg{A}$. There is a $\sigma \in
\Aut_{\rg{R}}/\ZM{n}$ acting like a cyclic permutation on $S = \{ x_1, x_2,
\ldots, x_t\}$.
\end{itemize}
Let $s \in \Z_{>1}$ and $\Phi_s(X) in \Z[ X ]$ be the \nth{s} cyclotomic
polynomial. If $\Psi(X) \in \ZM{n}[ X ]$ is a polynomial with $\Phi_s(X)
\equiv 0 \bmod (n, \Psi(X))$ and the extension $\rg{R} = \ZM{n}/(\Psi(X))$ is
simple Frobenius, we say that $(\rg{R}, \zeta, \sigma)$ is an \nth{s}
\textbf{cyclotomic extension} of $\ZM{n}$.

In general, if $\rg{R} \supset \ZM{n}$ is an algebra and $\zeta \in \rg{A}$
is such that $\overline{\Phi}_s(\zeta) = 0$, with $\overline{\Phi}_s(X) =
\Phi_s(X) \bmod n$, then we say that $\zeta$ is a \textbf{primitive} \nth{s}
root of unity modulo $n$.
\end{definition}
\begin{remark}
The reader may regard a cyclotomic extension $\rg{R}$ as an extension of the
ring $\ZM{n}$ which contains a primitive \nth{s} root of unity $\zeta$ and on
which an automorphism acts, that fixes $\ZM{n}$. One can prove - without
knowing that $n$ is prime - sufficient properties about $\rg{R}$ in order to
be allowed to work in the extension as if it was a finite field and $n$ were a
prime -- this behavior justifies the name of \textit{pseudo-fields} recently
employed by Lenstra.
\end{remark}
The pairs $(n,s)$ for which cyclotomic extensions exist are exceptional.  The
existence of such pairs is a strong property of $n$ with respect to $s$, that
often qualifies $n$ to behave like a prime. The following fact reflects this
claim: an \nth{s} cyclotomic extension of $\ZM{n}$ exists if and only if
\begin{eqnarray}
\label{fdiv}
  r \in \langle n \bmod s \rangle
  \quad \hbox{ for all }
  \quad r \mid n
\enspace.
\end{eqnarray}

Let $p$ be an odd prime and $k(p) = v_{p}\left(n^{p-1}-1\right)$, with
$v_{p}$ the $p$-adic valuation. If it exists, a \nth{p} cyclotomic
extension of $\ZM{n}$ may contain also a \nth{p^{k(p)}} primitive root
of unity; this is in fact true if $n$ is a prime. This leads to the
following
\begin{definition}
\label{dsat}
Let $p$ be a prime. The {\em saturation exponent} of $p$ is:
\begin{eqnarray}
\addtocounter{theorem}{1}
\label{3.53}
        k(p) & = & \left\{ \begin{array}{ll} 
        v_2(n^{2}-1)
        & \hbox{if $p = 2$ and $n = -1 \bmod 4$ } \\
        & \\
        v_{p}\left(n^{p-1}-1\right)
        & \hbox{ otherwise }
        \end{array} \right.
\end{eqnarray}

Let $m = \prod_{i} p_{i}^{e(i)} \in \N$ be the prime factorization of
an integer. The \emph{($n$-)saturated order above $m$} is:
\[
   \overline m = \prod_{i} \ p_{i}^{\max(e(i),k(p_{i}))} \enspace .
\] 
An \nth{m} cyclotomic extension $(\rg{R}, \sigma, \zeta)$ is called
\emph{saturated} if $m \geq \overline m$ and \emph{subsaturated} otherwise.
If $e(i) = k(p_{i})$ for all $p_{i} \mid m$, the extension is {\em
  minimal saturated}.
\end{definition}

Saturated extensions are characterized by the following property:
\begin{fact}
\label{tsat}
If $(\rg{R},\sigma,\zeta)$ is a saturated \nth{m} extension and $m' \mid
m^{h}$ for some $h > 0$ (i.e. $m'$ is built up from primes dividing
$m$), then $\rg{R}[X]/(X^{m'}-\zeta)$ is an \nth{m \cdot m'}
cyclotomic extension.

If $(\rg{R},\sigma,\zeta), (\rg{R'},\sigma',\zeta')$ are saturated \nth{m} and
\nth{m'} extensions for $(m, m') = 1$, then $(\rg{R} \times \rg{R}', \sigma
\circ \sigma', \zeta \cdot \zeta')$ is a saturated \nth{m m'} extension, for
the natural lifts of $\sigma, \sigma'$ to $\rg{R} \times \rg{R}'$.
\end{fact}

\cbstart
The use of saturated extensions in primality testing is given by the following 
\begin{lemma}[Cohen and Lenstra, \cite{CoLe}]
\label{lp}
Suppose that $p$ is a prime with $(p,n) = 1$, for which a saturated
\nth{p} cyclotomic extensions of $\ZM{n}$ exists. Then for any $r | n$
there is a $p$-adic integer $l_p(r)$ and, for $p>2$, a number $u_p(r)
\in \ZM{(p-1)}$, such that:
\begin{eqnarray}
\label{lpr}
r & = & n^{u_p(r)} \mod p \quad \mbox{and} \nonumber \\ & & \nonumber
  \\[-0.3cm] r^{p-1} & = & (n^{p-1})^{l_p(r)} \ \in \ \{1+p \cdot \Z_p
  \} \ \ \ \mbox{if} \quad p > 2, \\ & & \nonumber \\ [-0.3cm] r & = &
  n^{l_p(r)} \ \in \ \{1+2 \cdot \Z_2\} \quad \mbox{if} \quad p =
  2.\nonumber
\end{eqnarray}
\end{lemma}
\begin{proof}
Using \rf{fdiv}, the hypothesis implies that $r \in \ <n \mod p^k>$ for all 
$k \geq 1$ which implies \rf{lpr}.
\end{proof}
\cbend

Gauss and Jacobi sums over $\ZM{n}$ will be defined by means of characters
over saturated extensions. Let $p, q$ be two rational primes which do not
divide $n$, let $k > 0$ and $(\rg{R}, \zeta, \sigma)$ be a saturated \nth{p^k}
extension which additionally contains a primitive \nth{q} root of unity $\xi$;
the ring $\rg{R}$ need not be minimal with these properties. Let $\chi$ be a
multiplicative character $\chi: \ZMs{q} \rightarrow < \zeta >$ of conductor
$q$ and order $d | p^k$. If $d = 1$, $\chi$ is the trivial character $1$. The
(cyclotomic) Gauss sum of $\chi$ with respect to $\xi$ is
\[ \tau(\chi) = \sum_{x \in \ZMs{q}} \chi(x) \xi^x . \]
It can be shown that $\tau(\chi) \in \rg{R}^{\times}$, since $\tau(\chi) \cdot
\tau(\chi^{-1}) = \chi(-1) \cdot q$. For $a, b \in \Z$ such that $\chi^a,
\chi^b, \chi^{a+b} \neq 1$, the Jacobi sum
\[ j(\chi^a, \chi^b) = \sum_{x=2}^{q-1} \chi^a(x) \chi^b(1-x) = \frac{\tau(\chi^a) \cdot \tau(\chi^b)}{\tau(\chi^{a+b})} .\]

The {\sl multiple Jacobi-Sums $J_{\nu}(\chi)$ } are defined by:
\begin{eqnarray}
\label{2.31.1}
 J_1 & = & 1 \nonumber \\ J_{\nu + 1} & = &
	J_{\nu} \cdot j(\chi, \chi^{\nu}),\quad \hbox {for $\nu = 1,
	2, \ldots, d -2$} \\ J_{d} & = & \chi(-1) \cdot m \cdot
	J_{d-1} \nonumber
\end{eqnarray}

It is easy to verify by induction that:
\begin{eqnarray}
\label{2.31.2}
 J_{\nu} & = & \frac{\tau(\chi)^{\nu}}{ \tau(\chi ^ {\nu})}, \quad \hbox {for
	$\nu = 1, 2, \ldots, d$, where } \quad \chi^{d}=1.
\end{eqnarray}

Let $s = \prod_{ q \in Q} q$ be a product of primes from the set $Q$ such that
there is a $t = \prod_{p^k \in P} p^k$ with $P$ a set of prime powers and for
all $q \in Q, \ q-1 | t$.  Let $\rg{R}$ be the product of saturated \nth{p^k}
cyclotomic extensions and $C = \left\{ \chi_{p^k, q} \ : \ p^k \in P, q \in
Q\right\}$ be a set of characters of conductor $q$ and order $p^k$ with images
in $\rg{R}$. If $n = b(p^k) \cdot p^k + r(p^k)$ is the Euclidian division of
$n$ by each $p^k$, it can be shown \cite{BH, Mi1, Mi4} that a cyclotomic
\nth{s} extension of $\ZM{n}$ exists if
\begin{eqnarray}
\label{cpp}
  J_{p^k}^{b(p^k)}(\chi_{p^k, q}) \cdot J_{r(p^k)}(\chi_{p^k, q}) \in <
\zeta_{p^k} > , \quad \hbox{ for each } \quad \chi \in C.
\end{eqnarray}
Verifying these relations is the main stage of the CPP test. 
\begin{remark}
\label{cpp-comp}
Due to an analytic number theoretical Theorem of Pracher, Odlyzko and
Pomerance, one knows that two parameters $s, t$ can be chosen, such that $s >
\sqrt{n}$ and $t = O\left(\log(n)^{\log \log \log (n)}\right)$, while $s |
(n^t -1)$ for any $n$. The complexity of CPP is polynomial in $t$; both the
number of prime powers dividing $t$ and their size are upper bounded by
\begin{eqnarray}
\label{smallp}
B = O(\log \log(n)), \quad \omega(t) < B \quad \hbox{ and } \quad p^k || t
\Rightarrow p^k < B.
\end{eqnarray}
\end{remark}
We shall use an auxiliary construction involving \textit{dual elliptic primes}
in order to show that if $n$ passes the tests \rf{cpp} together with some
additional conditions - which are more involved to formulate, but can be
verified faster then \rf{cpp} - then either $n$ is prime, or it has a prime
factor $r$ with $l_{p^k}(r) = 1$ for all $p^k \in P$.

The constructions involve elliptic Gauss and Jacobi sums, which we shall
introduce below. We first define the simples analogue of cyclotomic extensions
for elliptic curves.

\begin{definition}
\label{el}
Let $n > 2$ be an integer and $\ell \nmid n$ be an odd prime. Let $\id{E}_n(A,
B)$ be an elliptic curve and $\psi_{\ell}(X)$ be the \nth{\ell} division
polynomial of the curve. Suppose that $F(X) \in \ZM{n}[ X ]$ is such that
\begin{itemize}
\item[1.] $F(X) | \psi_{\ell}(X)$.
\item[2.] If $\rg{E} = \ZM{n}[ X ]/(F(X))$ and $\Theta = X \bmod F(X) \in
\rg{R}$, then
\[ F( T ) = \prod_{i=1}^{(\ell-1)/2} (T - g_i(\Theta)) , \]
where $g_i(X)$ are the multiplication polynomials defined in \rf{mulpol}.  In
particular the elementary symmetric polynomials of $\Theta$ lay in $\ZM{n}$.
\end{itemize}

Then $F(X)$ is called an \textbf{Elkies factor} of $\psi_{\ell}(X)$ over
$\ZM{n}$ and $\rg{E}$ is an \textbf{Elkies ring}. Additionally, we let
\[ \rg{E'} = \rg{E}[ Y ]/\left(Y^2 - f(\Theta)\right) \quad \hbox{ and } 
\quad \Omega = Y \bmod (f(\Theta)
\in \rg{E'} \] be the \textit{two coordinates} Elkies ring.
\end{definition}
\cbstart 

Let $(\rg{R}, \zeta, \sigma)$ be a saturated \nth{\ell-1} cyclotomic extension
and $\chi: \ZMs{\ell} \rightarrow \rg{R}$ be a multiplicative character of odd
order.  We define Gauss sums in Elkies rings by:
\begin{eqnarray*}
\tau_e(\chi)  =  \sum_{i = 1}^{\ell-1} \chi(x) g_x(\Theta) 
\end{eqnarray*}
In the case when the order of $\chi$ is even and $\chi(-1) = -1$, the sums
above are vanishing due to the parity of $P_x$. One uses the $Y$ - coordinates
in the two coordinates Elkies ring, and some related multiplication
polynomials. The formal definition based on repeated addition of $P = (\Theta,
\Omega)$ in $\rg{E}'$ is in this case:
\begin{eqnarray*}
\tau'_e(\chi)  =  \sum_{i = 1}^{\ell-1} \chi(x) ([ x ] P)_{Y} 
\end{eqnarray*}
The values of $([ x ] P)_{Y}$ can be computed using $\Omega$ and polynomials
in $\Theta$; we skip the details here and refer to \cite{MMS, MV} for in depth
treatment of theoretical and computational aspects of elliptic Gauss and
Jacobi sums.

The Jacobi sums have no closed definition like in the cyclotomic case, so they
must be deduced as quotients of Gauss sums:
\begin{eqnarray*}
j_e(\chi^a, \chi^b) = \frac{\tau_e(\chi^a)\tau_e(\chi^b)}{\tau_e(\chi^{a+b})}
\quad \hbox{ iff } \tau_e(\chi^{a+b}) \in \rg{E}^{\times}.
\end{eqnarray*}

The case $\tau_e^{a+b}(\chi) \not \in \rg{E}^{\times}$ is improbable, but
cannot be excluded currently. This is best explained in the case when $n = r$
is a prime. Then $\id{E}_r(A, B)$ has a Deuring lift to some curve
$\id{E}_{\KH}(a, b)$. The Gauss sums of curves in characteristic $0$ have been
studied by R. Pinch in \cite{Pi} and it was shown that along with the ramified
primes dividing $\ell \cdot \Delta$, where $\Delta$ is the discriminant of the
curve, some spurious and unexplained primes may appear in the factorization of
the Gauss sum. Since $\Delta$ reduces to the discriminant of the curve
$\id{E}_r(A, B)$ which is non vanishing by definition and $\ell \neq r$, the
spurious primes may be divisors of $r$, in which case $\tau_e(\chi) \not \in
\rg{E}^{\times}$. If $n$ is not prime and $(\tau_e(\chi), n) \not \in \{1,
n\}$, a non trivial factor is found. We shall assume in our algorithm that the
case $(\tau_e(\chi), n) = n$ is scarce. It can be avoided by changing the
choice of $\ell$, as we shall detail below. If $\ell$ is a conductor, such
that $(\tau_e(\chi), n) = n$ for some character of conductor $\ell$, then we
say that $\ell$ is an {\em exceptional conductor} (for the curve $\id{E}_n(A,
B)$).

If $n = r$ is a prime, then $\Theta^r = g_{\lambda}(\Theta)$ for some $\lambda
\in \ZMs{\ell}$, an eigenvalue of the Frobenius. In that case, raising the
definition of the Gauss sum to the power $n$ yields:
\begin{eqnarray*}
\tau_e(\chi)^r & = & \left(\sum_{i = 1}^{\ell-1} \chi^r(x)
g_x(\Theta^r)\right) \\ & = & \sum_{i = 1}^{\ell-1} \chi^r(x) g_{\lambda x}
(\Theta) = \chi^{-r} (\lambda) \tau_e(\chi^r) \end{eqnarray*} and
\begin{eqnarray}
\label{potr}
\tau_e(\chi)^r/\tau_e(\chi^r)  =  \chi^{-r} (\lambda).
\end{eqnarray}
The right hand side of the equation can be computed, like in the cyclotomic
case by using multiple Jacobi sums in $\rg{R}$.  \cbend

\section{Elliptic extensions of rings}
\cbstart In this section we generalize the notion of \textit{cyclotomic
extension of rings} to elliptic curves.  We shall say that an Elkies algebra is
elliptic extension of $\ZM{n}$, if the power $n$ acts like a Frobenius,
i.e. \rf{potr} is verified when the prime $r$ is replaced by $n$. Note that
this is a slightly milder condition then the one for cyclotomic extensions,
since we are not interested in finding an actual factor of $F(X)$ which has
degree equal to the order of $n$ in the group $\ZMs{\ell}/\{-1,1\}$, i.e. the
degree of an irreducible factor of $F(X)$ in the case when $n$ is prime.
\begin{definition}
\label{primellext}
Let $m \in \N_{>2}$ be an elliptic Atkin pseudoprime: there is a curve
$\id{E}_m(A, B) : Y^2 = X^3 + A \cdot X + B$ associated to an order $\id{O}
\subset \K = \Q(\sqrt{-d})$ and such that $m = \mu \cdot \overline \mu$ for a
$\mu \in \id{O}$. Let $\ell$ be a rational prime such that
$\lchooses{-d}{\ell} = 1$:
\begin{itemize}
\item[A.] For each prime power $q || (\ell-1)$, there is a saturated \nth{q}
 cyclotomic extension $\rg{R}_q \supset \ZM{m}$. The rings $\rg{R}-q$ will
 also be called {\em working extensions}.
\item[B.] There is an \nth{\ell} cyclotomic extension $\rg{R}_{\ell} \supset
 \ZM{m}$ constructed by verifying \rf{cpp} over the extensions $\rg{R}_q$.
\item[C.] In particular, then 
\begin{eqnarray}
\label{lprq}
r \equiv n^{l_p(r)} \bmod \ell, \quad \hbox{ for $p | q$ a prime and for all }
\quad r | m.
\end{eqnarray}
\end{itemize}
Let $\psi_{\ell}(X)$ be the \nth{\ell} division polynomial associated to
$\id{E}_m(A, B)$ and suppose that an Elkies factor $F(X) | \psi_{\ell}(X)
\bmod m$ is known and $(\rg{E'}, \Theta, \Omega)$ is the two coordinates
Elkies algebra. For a prime power $q || (\ell-1)/2$ we let $\chi_q :
\ZMs{\ell} \rightarrow \rg{R}$ be a character of order $q$ and conductor
$\ell$. Suppose that:
\begin{itemize}
\item[1.] For each odd $q$, $(\tau_e(\chi), n) = 1$ and
\begin{eqnarray}
\label{oddtaue}
 \tau_e(\chi)^n/\tau_e(\chi^n) = \eta_q^{-n} \in < \zeta > .
\end{eqnarray}
\item[2.] For even $q$, $(\tau'_e(\chi), n) = 1$ and
\begin{eqnarray}
\label{eventaue}
 \tau'_e(\chi)^n/\tau'_e(\chi^n) = {\eta'}_q^{-n} \in < \zeta > . 
\end{eqnarray}

If the above conditions are met, we say that an \textbf{\nth{\ell}} elliptic
extension of $\ZM{n}$ related to $\rg{R}$ exists. The conditions
$\chi_q(\lambda) = \eta_q$ for odd $q$ and $\chi_q(\lambda) = \eta'_q$ for
even $q$ uniquely determine $\lambda_m \in \ZMs{\ell}$. This value will be
denoted as {\em the eigenvalue of the elliptic extension $\rg{E}$}.
\end{itemize}
\end{definition}
The point C. of the definition is a fact following from points
A. and B. and not a condition. The main fact about elliptic extensions is the
following:
\begin{theorem}
Let $n \in \N_{>2}$ be an integer and $\ell$ a prime not dividing $n$. If all
the conditions for existence of an \nth{\ell} elliptic extension of $\ZM{n}$
are fulfilled and $r | n$ is a prime, $\id{E}_r(\overline A, \overline B) =
\id{E}_n(A, B) \bmod r$, then
\begin{itemize}
\item[A.] The curve $\id{E}_r(\overline A, \overline B)$ has CM in $\id{O}$
and $\overline{F}(X) = F(X) \bmod r$ is an Elkies factor of its \nth{\ell}
torsion polynomial.
\item[B.] There is an eigenvalue $\lambda_r \in \F_{\ell}^{\times}$ of the
Frobenius of $\id{E}_r(\overline A, \overline B)$ such that $P^r = [ \lambda_r
] P$ for all points $P \in \id{E}_r(\overline A, \overline B)[ \ell ]$ such
that $\overline{F}(P_x) = 0$.
\item[C.] If $\lambda_m$ is the eigenvalue of the Elkies extension, $p$ is the
  prime dividing $ q$ and $l_p(r)$ is defined by \rf{lpr} with respect to the
  extension $\rg{R}$, then
\begin{eqnarray}
\label{ellext}
\chi_q(\lambda_r) & = & \chi_q(\lambda_m)^{l_p(r)} \quad \forall q. \\
\label{ellext2}
\chi_q(r/\lambda_r) & = & \chi_q(m/\lambda_m)^{l_p(r)} \quad \forall q. 
\end{eqnarray}
\end{itemize}
\end{theorem}
\begin{proof}
The point A. follows from Lemma 1. Point B. follows from the factorization
patterns of the division polynomial, e.g. \cite{Sch}, Theorems 6.1,
6.2. 

For proving \rf{ellext}, we use the fact that $\rg{R}$ is a cyclotomic
extension and let $\sigma : \zeta \mapsto \zeta^n$ act on the identities
\rf{oddtaue}. The $Y$ - component conditions \rf{eventaue} are treated
identically and will not be developed here.
\begin{eqnarray}
\label{bas}
\tau_e^{n^2}(\chi_q) & = & \left(\tau_e^n(\chi_q)\right)^n = \left(\eta_q^{-n}
\cdot \sigma(\tau_e(\chi_q)\right)^n \nonumber \\ & = & \eta_q^{-n^2} \cdot
\sigma(\eta_q^{-n} \tau_e(\chi_q^n)) = \eta_q^{-2n^2} \cdot
\sigma^2(\tau_e(\chi_q)), \ldots \\ \tau_e(\chi_q)^{n^k} & = & \eta_q^{-k n^k}
\sigma^k(\tau_e(\chi_q)).\nonumber
\end{eqnarray}
Inserting $k = \varphi(q)$ we obtain $ \tau_e^{n^{\varphi(q)}-1} = \eta_q $
and for $K = p \cdot \varphi(q)$, writing $N = n^{K}$, we have
\[ \tau_e(\chi_q)^{N-1} = 1. \]
If $r \mid n$ is a prime, by \rf{potr},
\begin{eqnarray*}
  \tau_e( \chi_q )^{r}\ \equiv \ \chi_q(\lambda_r)^{-r} \cdot (\tau_e(
		\chi_q^r ) ) \bmod r \cdot \rg{R}.
\end{eqnarray*}
Let $m \in \N$ be such that $m \equiv l_p(r) \mod p q$ and $m=u_{p}(r) \mod
 (p-1)$, with $u_{p} (r)$ and $l_p(r)$ defined by \rf{lpr}. Then $\sigma^m(
 \chi_q ) = \chi^r$ and
\begin{eqnarray}
\label{valus}
v_{\ell}(r-n^m)\ = \ v_{\ell} \bigl(n^m \cdot (r / n^m -1) \bigr) \geq
v_{\ell}(N-1).
\end{eqnarray}
We let $i = m$ in \rf{bas}, use $\sigma^m \bigl(\tau_e(\chi_q)\bigr) \ = \
\tau_e(\chi_q^r)$ and divide by \rf{valus}. This is allowed, since $(\tau_e(
\chi_q ), n) = 1$ by condition 1. Thus
\begin{eqnarray*}
	\tau_e( \chi_q )^{n^m-r} \ \equiv \ \bigl( \chi_q(\lambda_r) \cdot
		\eta_q^{-m} \bigr)^r \bmod r \cdot \rg{R}.
\end{eqnarray*} 
Raising this congruence to the power $a$, where $a$ is the largest divisor of
$(N-1)$ which is coprime to $\ell$, and using the above, we get :
\begin{eqnarray*}
1 \ \equiv \ \bigl( \chi_q(\lambda_r) \cdot \eta^{-m} \bigr)^{r \cdot a} \bmod
r \cdot \rg{R}.
\end{eqnarray*}
Since $(r a, \ell) = 1$, we deduce that $\chi_q( r ) \eta_q^{-m} \equiv 1
\bmod r \rg{R}$, and since $(\ell, n) = 1$ also $\chi_q( \lambda_r )
\eta_q^{-m} = 1$ and $\chi_q(\lambda_r) = \eta_q^m =
\eta_q^{l_p(r)}$. This holds for all primes $r \mid n$ and by
multiplicativity, for all $r | n$. In particular, since $ l_p ( m ) = 1$,
it follows that $\chi_q ( \lambda_m ) = \eta_q$, thus recovering the
definition of the eigenvalue of the elliptic extension. The proof of
\rf{ellext} is complete. As for \rf{ellext2}, it follows from \rf{ellext} and
\rf{lprq}. 
\end{proof}
\cbend

The notion of elliptic extension for composites is now straight forward:
\begin{definition}
\label{compellext}
Let $L = \prod_{i=1}^k \ell_i \in \N$ be square-free, with $\ell_i$ being
primes. Assume that there is an \nth{\varphi(L)} saturated working extension
$\rg{r}_L \supset \ZM{m}$ and an \nth{L} extension $\rg{R}_L \supset
\rg{r}_L$. 

Suppose also that $m$ is Atkin pseudoprime so there is a curve
$\id{E}_m(A, B)$ associated to an order $\id{O} \subset \K = \Q[ \sqrt{-d} ]$.
We say that an \nth{L} elliptic extension exists, if the conditions of
Definition \ref{primellext} are fulfilled for all $\ell_i$ in the working
extension $\rg{r}_L$ or subextensions thereof.
\end{definition}

Note that relation \rf{ellext} is a strengthening of the consequence
$p \equiv m^k \bmod L$, usual in classical cyclotomy tests. It follows from
the definition and \rf{ellext} \rf{ellext2} that 
\begin{eqnarray}
\label{compext}
\lambda_r & \equiv & \lambda_m^{k_L(r)} \bmod L \quad \hbox{ for $k_L(r) \equiv
  l_p(r)$, for all $ p | \varphi(L)$}, \\
\label{compext2}
(r/\lambda_r) & \equiv & (m/\lambda_m)^{k_L(r)} \bmod L.
\end{eqnarray}
We shall combine this strengthening with properties of dual elliptic
pseudoprimes, which we introduce in the next section, with the goal of
eliminating the final trial division \rf{findiv} in cyclotomy tests of a given
pair of dual elliptic primes.

\section{Dual Elliptic Primes and Pseudo-primes}
We start with the definition of the dual elliptic primes, which is, as
mentioned in the introduction, related to the notion of twin primes in the
rational integers.
\begin{definition}
We say that two primes $p$ and $q$ are {\sl dual elliptic primes}
associated to an order $\id{O} \subset \K = \Q(\sqrt{-d})$, if there
is a prime $\pi \in \id{O}$ such that $p = \pi \cdot \overline {\pi}$
and $q = (\pi + \varepsilon) \overline {(\pi + \varepsilon)}$ with
$\varepsilon = \pm 1$.
\end{definition}
 Dual elliptic primes \textit{exist}: In the ECPP program, a special flag was
   introduced in order to skip dual pseudoprimes, which do not reduce the size
   of the numbers to be proved prime; it happens regularly that the flag is
   set \cite{Mo4}. Furthermore, empirical considerations of Galbraith and McKee
   \cite{GaK} suggest they are sufficiently frequent, in order to develop
   efficient algorithms in which they are used. The problem of showing that
   dual elliptic primes have a satisfactorily asymptotic distribution is
   certainly much harder.

We define in the spirit of pseudoprimality followed from the introduction, a
pair of dual elliptic pseudoprimes as follows:

\begin{definition}
\label{DEPSP}
Let $m$ and $n$ be two strong pseudoprimes, $\id{O} \subset \K =
\Q(\sqrt{-d})$ an order in an imaginary quadratic extension and assume
that there are two curves $\id{E}_m(A, B), \ \id{E}_n(C, D)$ which are
both associated to $\id{O}$ in the sense of Definition \ref{Atkassoc}. In
particular, $m, n$ are Atkin pseudoprimes. Furthermore, we assume
that:
\begin{itemize}
\item[1.] There are a point $P \in \id{E}_m(A, B)[ n ]$ and a point $Q \in
  \id{E}_n(C, D)[ m ]$ and the (Atkin) - sizes of the curves are
\[  \left|\id{E}_m(A, B)\right| = n, \quad \hbox{ and } \quad
\left|\id{E}_n(C, D)\right| = m . \]
\item[2.]  The sizes $m$ and $n$ factor in $\id{O}$ as
\begin{eqnarray}
\label{dualfact} 
m = \mu \overline \mu, \quad \hbox{ and } \quad n = (\mu + \varepsilon) \cdot
\overline{(\mu + \varepsilon)}, \quad \hbox{ with } \quad \varepsilon = \pm 1.
\end{eqnarray}
Note that from \rf{coprim} we have that $m, n$ are square-free.
\item[3.] The polynomial $H_{\id{O}}(X)$ has a root $j_m$ modulo $m$, and a
  root $j_n$ modulo $n$, and the curves $\id{E}_m(A, B), \id{E}_n(C, D)$ have
  invariants which are rational functions in these values.
\item[4.]  Both $m$ and $n$ have no prime factor $p < 5$.
\end{itemize}
\vspace*{0.5cm} If these conditions are fulfilled, the pair $(m, n)$ is called
a pair of {\em dual elliptic pseudoprimes} associated to the order $\id{O}$.
\end{definition}

Finding a point $P$ on $\id{E}_m(A, B)$ can be done by adapting a trick of
\cite[8.6.3]{AtMo}, thereby bypassing the problematic extraction of a square
root modulo $m$. This works as follows: find $x_0\bmod m$ for which $\lambda =
x_0^3 + a x_0 + b\mod m$ is such that $\legendre{\lambda}{m}=1$. Then
$P=(\lambda x_0, \lambda^2)$ is a point on the curve $Y^2 = X^3+A \lambda^2
X+B \lambda^3$, which should be isomorphic to $\id{E}_m(A, B)$ if $m$ is
actually a prime\footnote{I thank F. Morain for this observation}.

Practically, dual elliptic pseudoprimes are found by featuring a pair of
strong pseudoprimes $(m,n)$; the pseudoprime test may consist in taking the
roots \linebreak $\sqrt{ - d \bmod m}, \sqrt{ - d \bmod n}$, operations which
are anyhow necessary in the context. The integers $m$ and $n$ both split in a
product of two principal primes in $\K$, such that there is a pair of factors
which differ by $\pm 1$. Once such pseudoprimes are found, the invariants
$j_{m}, j_{n}$ must be computed by methods explained in \cite{Mo},
\cite{AtMo}. Then the curves $\id{E}_m(A, B), \id{E}_n(C, D)$ can be built and
points on these curves are chosen as explained above. The points are used in
order to perform an elliptic pseudoprime test, as required in point 1 of the
Definition \ref{DEPSP}. In practice one notes that, given a strong pseudoprime
$n$, finding an appropriate order $\id{O}$ and a dual elliptic pseudoprime $m$
to $n$ is a particular form of the first round of an elliptic curve primality
test (ECPP) \cite{Mo}. In particular, the heuristic arguments based upon
\cite{GaK} suggest that this step requires cubic time.

The easiest fact about dual elliptic pseudoprimes is the following:
\begin{lemma}
\label{2.2}
Two dual elliptic pseudoprimes $(m,n)$ associated to an order $\id{O}$ are
simultaneously prime or composite. Furthermore, if $m, n$ are composite and
$\id{O} \subset \K = \Q(\sqrt{-d})$, then for any prime divisor $\ell | m
\cdot n$ there is a $\lambda \in \id{O}(\K)$ such that $\ell = \lambda \cdot
\overline \lambda$.
\end{lemma}
\begin{proof}
Assume $m$ is prime. Then item 1. of the Definition \ref{DEPSP} requires also
an elliptic Fermat primality proof for $n$. It implies that for any possible
prime $q | n$, the curve $\id{E}_q(\overline A, \overline B) = \id{E}_n(A, B)
\bmod q$ has a point of prime order $m > (\sqrt{n} - 1)^2$. This cannot hold
for primes $q < \sqrt{n}$ and thus $n$ is prime too. Conversely, if $n$ is
prime, $m$ is also prime by the same argument. This confirms the first
statement.

Suppose now that $m$ and $n$ are composite and $\ell \in \N$ is a prime so
that $\ell | n$, say. The condition \rf{coprim} implies that $n$ is square -
free and Lemma \ref{l1} together with point 2. of the Definition
\ref{Atkassoc} imply that $\ell$ splits in a product of principal ideals of
$\id{O}$, which completes the proof.
\end{proof}

We shall assume from now on, without restriction of generality, that
$\varepsilon = 1$ in the Definition \ref{DEPSP} (note that changing
the sign of $\varepsilon$ amounts to interchanging $m$ and $n$). We
prove that the tests required by the definition imply that, if
dual elliptic pseudoprimes are composite, then their least prime
factor {\em has } the dual elliptic prime property.

\begin{theorem}
\label{sfac}
 Let $(m,n)$ be a pair of composite dual elliptic pseudoprimes associated to
 an order $\id{O} \subset \Q(\sqrt{-d})$ and let $p \mid m$ be the least prime
 factor of $m$. Then there is a prime factor $q \mid n$, such that $p, q$ are
 dual elliptic primes.

Furthermore, if the prime $q$ is not the least prime factor of $n$,
 then both $m$ and $n$ are built up of at least three prime factors.
\end{theorem}
\begin{proof}
By Definitions \ref{Atkps} and \ref{DEPSP}, there is a point $P \in
\id{E}_m(A, B)$ with $[ n ] P = \eu{O}$. Let $\rg{P} = P \bmod p \in
\id{E}_p(\overline A, \overline B) = \id{E}_m(A, B) \bmod p$; it has an order
$h \mid n$. If $h$ is a prime, then $p, h$ are dual elliptic primes and the
proof is completed. Let us thus assume that $h$ is composite and $q \mid h
\mid n$ is the least prime dividing $h$, so $h = q \cdot u$, with some $u >
1$. By the choice of $q$ it follows that $q^2 < qu = h$. We then consider $Q
\in \id{E}_n(C, D)[ m ]$ and the point 
\[ \rg{Q} = \left( Q \bmod q \right) \in
\id{E}_q(\overline C,\overline D) = \left(\id{E}_n(C, D) \bmod q\right), \]
which must have a non trivial order $h' \mid m$, since $Q$ is an \nth{m}
torsion point. The choice of $p$ implies $h' \geq p$. Applying the Hasse
inequalities to $h$ and $h'$ we find:
\begin{eqnarray*}
\begin{array}{r c c c l}
q^2 & \leq & q \cdot u & \leq & (\sqrt{p}+1)^2, \\
(\sqrt{p}-1)^2 & \leq & h & \leq & u^2,  \\
p & \leq & h' & \leq & (\sqrt{q}+1)^2.
\end{array}
\end{eqnarray*}
Thus, from the first two lines, $q \leq \sqrt{p}+1 \leq u +
2$ and combining to the other inequalities we have:
\[ q \cdot u \leq (\sqrt{p}+1)^2 \leq (\sqrt{q}+1)^2 + 1 + 2\sqrt{p}. \]
After division by $q$, we find the following bonds on $u$:
\[ 1 \leq u \leq \frac{(\sqrt{q}+1)^2}{q}+\frac{2u+3}{q} < 1 + 4/q + 2/\sqrt{q}
+ 2u/q, \] and since $q \geq 5$, also $3 u/5 < 3$. This is impossible, since
$u > q \geq 5$ is an integer. Thus $u = 1$ and $h = q$ is prime, which
completes the proof of the second statement.

We had chosen $q$ as the least prime factor of $h$, the order of the
point $\rg{P} \in \id{E}_p(A, B)$. We now show that if $q$ is not the least
prime factor of $n$, then $n$ has more then two prime factors. Assume
that $q' < q$ is the least prime dividing $n$. By the proof above,
there is a prime $p'\mid m$ such $(q', p')$ are dual; also the
premises imply that $p'> p$. Given the double duality, we have the
following factorizations in $\id{O}(\K)$:
\begin{eqnarray*}
\begin{array}{r c l c r c l c l}
p & = & \pi \cdot \overline \pi & ; & q & = & \rho \cdot \overline \rho & = &
(\pi + \delta)\left(\overline{\pi + \delta}\right) \\ p' & = & \pi' \cdot
\overline \pi' & ; & q' & = & \rho' \cdot \overline \rho' & = & (\pi' +
\delta')\left(\overline{\pi' + \delta'}\right),
\end{array}
\end{eqnarray*}
where $\delta, \ \delta' = \pm 1$ and $\pi$ and $\pi'$ can be chosen such that
their traces be positive.

We assume that $m = p \cdot p'$ and $n = q \cdot q'$ and insert the
last equations in the factorizations of $m$ and $n$ in $\K$:
\begin{eqnarray*}
\begin{array}{c c c c c c c}
m & = & \mu \cdot \overline \mu & \hbox{and} & \mu & = & \pi \cdot
\pi'\\ n & = & (\mu+1) \cdot (\overline \mu + 1) & \hbox{and} & \mu + 1 &
= & (\pi+\delta) \cdot (\pi'+\delta' ).
\end{array}
\end{eqnarray*}
Subtracting the right hand side equations, we find $1 - \delta \cdot \delta' =
\delta \pi' + \delta' \pi$. If $\delta = \delta'$, this implies $\pi + \pi' =
0$ and $\mu$ is a square. If $\delta = -\delta'$ then $\pi' - \pi = 2 \delta$
and $\rho = \pi+\delta$, so $\rho' = \pi' + \delta' = \pi + 2 \delta + \delta'
= \pi + \delta = \rho$, then $\nu$ is a square. But both $\mu, \nu$ were
assumed square-free, a contradiction which confirms that at least one of $m$
and $n$ must have three factors.

Assume now that one of $m, n$ is built up of two primes, say $m = p \cdot p'$,
while $n = q \cdot q' \cdot q''$, where $q''$ is a factor which may be
composite and $q' < q < q''; p < p'$. By duality, we have $q' >
(\sqrt{p'}-1)^2$ and $q'' > q > (\sqrt{p}-1)^2$, thus 
\[ n =  q \cdot q' \cdot q'' > m \cdot \left((p+1 - 2\sqrt{p}) \cdot
(1-2/\sqrt{p})(1-2/\sqrt{p'})\right) .\]
For $p' > p \geq 11$ it follows that $n > 1.367 \ m$ and $m > 121$, in
contradiction with $n < m+1 + 2/\sqrt{m} < 1.2 \ m$. The remaining cases can be
eliminated individually, using the fact that small primes $5 \leq p < 11$ split
in principal ideals only in few imaginary quadratic extensions, and in those
cases, if $p = \pi \cdot \overline \pi$, then $\pi \pm 1$ is not prime.
\end{proof}
An immediate consequence is the following:
\begin{corollary}
\label{cdif}
Let $(m,n)$ be dual elliptic pseudoprimes associated to the order $\id{O}
\subset \K = \Q(\sqrt{-d})$ and $k = k(m,n) = \max\{\Omega(m), \Omega(n)\}$,
where $\Omega(x)$ denotes the number of prime factors of $x$, with
repetition. Then there are two primes $p \mid m$ and $q \mid n$ such that:
\begin{eqnarray}
\label{2.4}
 | p - q | < 2 \sqrt{\max(p,q)} < 2 {\root 2k \of {\max(m,n)}} \leq 2 {\root 4
   \of {\max(m,n)}} .
\end{eqnarray}
\end{corollary}
\begin{proof}
Suppose that $m$ has $k=k(m,n)$ factors and let $p$ be its least prime
factor, so $p < m^{1/k}$. Let $q$ be the dual prime of $p$
dividing $n$: the existence of $q$ follows from the previous
theorem. Then \rf{2.4} follows from the duality of $p$ and $q$ and
the bound on $p$.
\end{proof}
We finally show that dual elliptic primes with two factors might
exist. This leads to a formula which reminds formulae for the prime
factors of Carmichael numbers.
\begin{theorem}
\label{twofac}
 Let $(m,n)$ be a pair of dual elliptic pseudoprimes associated to an order
 $\id{O} \subset \Q[\sqrt{-d}]$ and suppose that both are built up of exactly
 two prime factors. Let $m = \mu \cdot \overline \mu$ and $n = (\mu+1) \cdot
 (\overline \mu + 1)$ be the factorizations of $m$ and $n$ in $\K$. Then there
 is a prime $\pi$, an element $\alpha \in \id{O}$ and a unit $\delta$,
 such that:
\begin{eqnarray}
\label{2.5}
\nu & = & (\pi + \delta) \cdot (\alpha \pi + \delta) \quad \hbox{and} \\
\nonumber \mu & = & \pi \cdot (\alpha(\pi + \delta) + \delta) .
\end{eqnarray}
\end{theorem}
\begin{proof}
  Let $m = p \cdot p'$ and $n = q \cdot q'$ be the rational prime
  factorization of $m$ and $n$. Since $m$ and $n$ have only two prime factors,
  it follows from Theorem \ref{sfac} that the least primes, say $p, q$ must be
  dual to each other. So let $p = \pi \cdot \overline \pi$ and $q = \rho \cdot
  \overline \rho = (\pi+\delta) \cdot (\overline{\pi + \delta})$.
  
  Let also $p' = \pi' \cdot \overline \pi'$ and $q'=\rho' \cdot \overline
  \rho'$. The size of $\id{E}_{q'}(\overline C, \overline D) = \id{E}_n(C, D)
  \bmod q'$ divides $m$ and it follows, after an adequate rearrangement of
  conjugates, that there is an $\varepsilon = \pm 1$ such that $\rho' +
  \varepsilon$ is divisible by either $\pi$ or $\pi'$. 

If the divisor was $\pi'$ we would reach a contradiction like in the last step
  of the proof of Theorem \ref{sfac}. Assume thus that $\rho' = \alpha \pi -
  \varepsilon$, the divisor being $\pi$. Symmetrically, $\pi' = \beta \rho +
  \varepsilon'$. First consider the splitting of $\nu$:
\[ \mu + 1 = \nu = \rho \cdot \rho' = (\pi + \delta)(\alpha \pi + \varepsilon) 
= \pi (\alpha \pi + \alpha \delta + \varepsilon) + \varepsilon \delta
\] Reducing the above equation modulo $\pi$, we conclude that
$\varepsilon \delta = 1$ and thus $\varepsilon = \delta$, both factors
being $\pm 1$. Let us compare the two expressions for $\mu$:
\[ \mu = (\alpha \pi^2 + \delta(\alpha+1) \pi + 1) - 1 =
\pi (\beta(\pi+\delta) + \varepsilon') \] and, after dividing
$\pi$ out,
\[  (\alpha - \beta) (\pi + \delta) = \varepsilon' - \delta. \]
If $\varepsilon' = \delta$, then $\alpha = \beta$ and the claim follows.
If $\alpha \neq \beta$, one can divide both sides by $\alpha - \beta$:
\[      \pi+\delta = \pm \frac{2}{\alpha - \beta }, \quad \hbox{ thus } \quad
(\alpha - \beta) | (2). \]
Assuming that $\alpha - \beta = \zeta \in \id{O}(\K)^{\times}$, one finds $\rho
= \pi + \delta = 2 \zeta'$, for some related root of unity $\zeta'$. This
contradicts the fact that $\rho \overline \rho = q \geq 5$.

Finally we have to consider the case when $\alpha - \beta \in \id{O}$ divides
$2$ and is not a unit. The only quadratic imaginary extension in which the
prime $2$ factors in principal ideals is $\K = \Q[ \imath ]$. Thus for $\K
\neq \Q[ \imath ]$ we must have $\alpha = \beta$ and the statement
follows. Finally, if $\K = \Q[ \imath ]$, we substitute $\alpha - \beta = 1
\pm \imath$ in the previous identity and find solutions for $\pi, \pi';
\rho, \rho'$ which are also of the shape \rf{2.5}; this completes the proof.
\end{proof}

\subsection{Elliptic extentions of dual elliptic pseudoprimes}
Let $(m, n)$ be a pair of dual elliptic pseudoprimes associated to an order
$\id{O} \subset \K = \Q(\sqrt{-d})$ and $\id{E}_m(A, B), \id{E}_n(C, D)$ be
the respective curves. We have shown that to the least prime $p | m$ there is
a dual elliptic $q | n$ and both factor into principal primes in
$\Q(\sqrt{-d})$; let $p = \pi \cdot \overline \pi$ and $(\pi +
\delta)(\overline \pi + \delta) = q$ be these factorizations, with $\delta =
\pm 1$. Suppose that $L$ is a square free integer for which the $L$ - torsions
of the curves $\id{E}_m(A, B)$ and $\id{E}_n(C, D)$ give raise to elliptic
extensions of $\ZM{m}, \ZM{n}$. Let these extensions be defined over the
saturated \nth{\varphi(L)} cyclotomic extensions $(\rg{R}_m, \zeta_m,
\sigma_m)$ and $(\rg{R}_n, \zeta_n, \sigma_n)$ respectively.

If $m, n$ are primes, then the eigenvalues of the Frobenius are $\mu + 1,
\overline \mu + 1$ for $\Phi_m$ and $\mu, \overline \mu$ for $\Phi_n$, as one
deduces from the sizes of the curves. By definition of the Elkies primes, they
split in $\id{O}(\K)$ and for each prime $\ell | L$ we have $(\ell) = \id{L}_1
\cdot \id{L}_2$; one should check additionally that:
\begin{eqnarray}
\label{ev}
\lambda_m & \in &  \left\{ \mu + 1 \bmod \id{L}_1, \overline \mu + 1 
\bmod \id{L}_1\right\} , \\
\lambda_n & \in &  \left\{ \mu \bmod \id{L}_1, \overline \mu 
\bmod \id{L}_1\right\} . \nonumber.
\end{eqnarray}
Then \rf{compext} implies that there are two integers $k, k'$ such that:
\[ \pi \equiv \mu^k \bmod L \id{O} \quad \pi + \delta \equiv (\mu+1)^{k'} 
\bmod L \id{O}. \] 
\begin{remark}
\label{detk}
The numbers $k, k'$ are determined by $k \equiv l_{v^i}(p)$ and
$k' \equiv l_{v^i}(q)$ for each prime power $v^i || \varphi(L)$. Using also
\rf{lprq} both for $m$ and $n$, it follows that
\begin{eqnarray}
\label{ellk}
    (\mu+1)^{k'} - \mu^k \equiv \delta \bmod L \id{O} . 
\end{eqnarray}
Note that the fact that the \nth{\varphi(L)} extension is saturated requires
in particular, that for each prime $v | \varphi(L)$ with saturation exponent
$j$, the power $v^j | \varphi(L)$.
\end{remark}
One may consider \rf{ellk} as an equation in the unknowns $k, k'$. In
particular, $(1, 1)$ is always a possible solution, for which $\delta = 1$.
It is possible that for certain $L$, the trivial is the only solution. We
shall say that a square free integer $L$, which is product of primes $\ell$
which split in $\id{O}(\K)$ and such that \rf{ellk} has only the trivial
solution is a \textbf{good $L$} -- with respect to the dual pseudoprimes $m,
n$. This property has important consequences for the cyclotomy test as shown
by the following
\begin{theorem}
\label{thL}
Let $m, n$ be dual elliptic pseudoprimes associated to an order $\id{O}
\subset \K = \Q[\sqrt{-d}]$ and let $m = \mu \cdot \overline \mu, \ n = (\mu +
1) (\overline \mu + 1)$ be the respective factorizations in $\id{O}$.

Suppose that $L \in \N$ is a square free integer for which an \nth{L} elliptic
extension exists both for $\ZM{m}$ and $\ZM{n}$ and they are defined using the
saturated $\varphi(L)$ extensions $(\rg{R}_m, \zeta_m, \sigma_m)$ and
$(\rg{R}_n, \zeta_n, \sigma_n)$ respectively; suppose that \rf{ev} holds for
the eigenvalues of these extensions. If the system \rf{ellk} has only the
trivial solution $(k, k') = (1, 1)$ and $p \ | \ m;  q \ | \ n$ are two dual
elliptic primes, then
\begin{eqnarray}
\label{not}
l_v(p) \equiv l_v(q) \equiv 1 \bmod v^N, \quad \hbox{ for each prime $v |
  \varphi(L) $ and $N > 0$}.
\end{eqnarray}
\end{theorem}
\begin{proof}
The statement \rf{not} is a direct consequence of Remark \ref{detk} and the
fact that the \nth{\varphi(L)} extensions is saturated. 
\end{proof}

The Theorem suggests the following procedure for eliminating the final trial
division step in the cyclotomy test:
\begin{itemize}
\item[1.] Start with a pair of dual elliptic pseudoprimes $m, n$ associated to
  an order $\id{O}$ and choose two parameters $s, t$ with $s | (n^t-1, m^t-1)$
  for a cyclotomy test, as indicated by Remark \ref{cpp-comp}.
\item[2.] Search by trial and error a square free $L$ such that $t |
  \varphi(L)$ and an elliptic \nth{L} extension of $\ZM{m}$ and of $\ZM{n}$
  exists. Note that the primes dividing $L$ need to be Elkies primes, which
  depends on $\id{O}$ and not on the individual values of $m, n$. They may but
  need not divide $s$.
\item[3.] Suppose that additionally \rf{ev} holds and \rf{ellk} has only the
  trivial solution.
\end{itemize}
If such a construction succeeds together with the main stage of the cyclotomy
test for $m, n$ and these are not primes, then there are two (dual elliptic)
primes $p \ | \ m; q \ | \ n$ with $p < \sqrt{m}, q < (\sqrt{p}+1)^2$ and such
that
\begin{eqnarray}
\label{trdivce}
p \equiv m \bmod L \cdot s \quad \hbox{ and } \quad q \equiv n \bmod L \cdot s.
\end{eqnarray}
This follows from \rf{not} together with the fact that the existence of an
\nth{Ls} cyclotomic is jointly proved by the cyclotomy test and the above
additional steps. In particular, the final trial division is herewith
superfluous.
\subsection{Heuristics}
We complete this section with a heuristic analysis for the odds of finding $L$
which verifies the conditions of Theorem \ref{thL}. We start with some
simplifications and consider one prime $\ell | L$ with $\ell > 3$ and which
factors in $\id{O}$ according to $(\ell) = \id{L}_1 \cdot \id{L}_2$. We let $x
\equiv \mu \bmod \id{L}_1$ and $y \equiv \overline \mu \bmod \id{L}_1$, with
$x, y \in \F_{\ell}^{\times}$. Restricted to $L = \ell$, the system \rf{ellk}
becomes in this notation: $x^k + \delta = (x + 1)^{k'}$ and $y^k + \delta = (y
+ 1)^{k'}$. Fix a generator $g \in \F_{\ell}^{\times}$ and consider the
discrete logarithm in $\F_{\ell}^{\times}$ with respect to $g$. 

We shall assume for simplicity that $x, y, x+1, y+1$ also generate the
multiplicative group $\F_{\ell}^{\times}$, so
\begin{eqnarray}
\label{gens}
\log(a) \in \ZMs{\ell-1} \quad \hbox{ for } \quad a \in \{x, y, x+1, y+1\}.
\end{eqnarray}
Consider the functions $f_x, f_y : \ZM{\ell-1} \rightarrow \ZM{\ell-1}$ given
by
\[ f_x(k) = \frac{\log(x^k + \delta)}{\log(x+1)} \quad   f_y(k) =
\frac{\log(y^k + \delta)}{\log(y+1)} .\] The system \rf{ellk} is now $f_x(k) =
f_y(k) = k'$. We exclude the couple $(1, 1)$, corresponding to the trivial
solution, from the graph of $f_x$. Furthermore $x \neq 0 and x + 1
\neq 0$, and thus $x^k \neq -\delta$ and $(x+1)^{k'} \neq \delta$. This
excludes an additional pair $(a,b)$ from the graph of $f_x$. The same holds
for $f_y$ and both maps are restricted to domains and codomains of equal size
$\ell-3$.

\begin{fact}
\label{heur}
Our heuristic is based on the assumption that the functions $f_x, f_y$ are
well modeled by random permutations of $S_{\ell-3}$. In particular, modulo a
redefinition of either domain or codomain, the maps are invertible and the
system \rf{ellk} reduces to $f^{-1}_y \circ f_x(k) = k$. According to our
model, the map $h_{x, y} = f^{-1}_y \circ f_x$ is also a random permutation
and it should have at least one fixed point. 

The number of fixed points of random permutations is well understood: it has
expected value $1$ and is Poisson distributed. Asymptotically, the individual
probabilities $P_k = P(h \hbox{ has $k$ fixpoints}) \rightarrow \frac{1}{k
!}$. Along with the expected value, we are interested in the probability that
$h$ has no fix points at all, which is $P_0 = 1/e$. The $\limsup_{X
\rightarrow \infty} \frac{X }{\varphi(X) \log \log(X)} \leq C$ for some $C >
0$, \cite{}. For fixed $x, y, x+1, y+1$ and a given $0 < B < \log(m)$, a prime
$\ell \equiv 1 \bmod B$ such that \rf{gens} holds, occurs with probability $P
> C'/(\log\log(B))^4$. The heuristic model implies that with expectation
$1/e$, \rf{ellk} will only have the trivial solution for such a prime.

A further approach which can be analyzed with the same model is the following:
choose $\ell_1, \ell_2$ like above, and let $n_1, n_2$ be the respective
number of fixed points. The expected values are $n_1 = n_2 = 1$. Suppose that
$(\ell_1 -1, \ell_2 -1) = d$ and let $k, k' \in \ZM{\varphi(L)}$ be a non
trivial solution of \rf{ellk}. Let $k_i \equiv k \bmod \ell_i-1; \ k'_i \equiv
k' \bmod \ell_i-1, \ i = 1, 2$ be the exponents with respect to $\ell_i$. They
correspond to some of the $n_i$ solutions modulo $\ell_i$, and thus $k_1
\equiv k_2 \bmod d; \ k'_1 \equiv k'_2 \bmod d$. Since there is in average
only one solution modulo each prime, this solution must verify the above pair
of additional conditions, which are met with probability $1/d^2$. Thus, if $d
> 1$, the probability that \rf{ellk} has a solution for $L$ as above is $1/d^2
< 1/e$ and trying at least two primes yields a stronger filtering.
\end{fact}

Certainly, the condition \rf{gens} is only necessary for a simpler heuristic
argument. The analysis may become difficult when some of $x, y, x+1, y+1$ are
not generators. The odds of finding a \textit{good} $L$ are though the same
range of magnitude. For the purpose of finding good $L$, we thus propose the
more general algorithm:
\begin{quotation}
\NI {\bf Algorithm ACE}{\em ( Auxiliary Cyclotomic \& Elliptic Extensions ) }
\vspace*{0.3cm}
\\
\textit{Input}. $m, n$ a couple of dual elliptic pseudoprimes with respect
to $\id{O}$ with given factorization; $t$, an exponent for a CPP test.
\textit{Output} $L$ a square-free integer with $t | \varphi(L)$ and such that
\rf{ellk} has only the trivial solution modulo $L$.
\vspace*{0.1cm}
Compute a sequence of primes $\ell_i > 3; \ i = 1, 2, \ldots h$ and
  let $L_i = \prod_{j \leq i} \ell_i$, such that
\begin{itemize}
\item[(i)] $d_i = (\ell_i, \varphi(L_{i-1})) > 1$.
\item[(ii)] $L = L_h$ is such that $t | \varphi(L)$.
\item[(iii)] The equations \rf{ellk} have no non trivial solutions modulo $L$.
\end{itemize}
  \vspace*{0.1cm}
\end{quotation}

\begin{remark}
\begin{itemize}

\item[A.] We have implemented this algorithm. In most cases, the equations
\rf{ellk} had only the trivial solution for $L$ a product of two primes. In
more then one fourth of the cases, this happened already for one prime, and we
encountered no case in which a product of more then three primes was necessary
for a good $L$. Thus the experimental results in the general case are close to
the heuristic predictions for the particular case in which \rf{gens} holds.
\item[B.] The condition (i) has the following purpose: in general, we reach a
  good $L_j$ already for $j \leq 3$, however the condition $t | \varphi(L)$
  will not be fulfilled. Suppose thus that $L_j$ is good and \rf{ellk} has at
  least one non trivial solution $(k, k')$ for $\ell_{j+1}$. If $d_{j+1} > 1$,
  since $L_j$ is good, we must have $k \equiv k' \equiv 1 \bmod d_{j+1}$: this
  allows filtering. In practice, one shall choose $d_{j}$ to be at least
  divisible by some factors of $t$.

\item[C.] Assume that $B > 0$ is such that all prime power factors of $t$ are
  $< B$ and the number of prime factors is also $< B$ -- see \rf{smallp}. We
  claim that the Algorithm ACE will complete in average time
  $O(B^{3+\varepsilon})$. For the analysis, we use again the slower approach,
  in which one seeks for each prime power $v | t$ a good prime $\ell \equiv
  \bmod \ v$, such that \rf{gens} holds. By Fact \ref{heur}, a good prime for
  which \rf{gens} holds occurs with probability $O(1/\log
  \log^4(B))$. Combining with the probability to find a prime $\ell \equiv 1
  \bmod v$ estimated with the Linnick constant, we deduce that for
  sufficiently large $t$ and thus $B$, there is a good prime $\ell(v) <
  B^{1+\varepsilon}$ with $\ell(v) \equiv 1 \bmod v$ for each prime power $v |
  t$. It tales $O(B^{2+\varepsilon})$ to find such a prime. Repeating this for
  all $v | t$ will take at most $O(B^{3 + \varepsilon})$ operations, as
  claimed. In practice, by \rf{smallp}, $B = O(\log \log(m))$ and thus the
  time required by the ACE Algorithm is negligeable.
\end{itemize}
\end{remark}
\subsection{On constructing Elkies factors}
We finally add some detail on the construction of the Elkies factors of
\nth{\ell} torsion polynomials $\psi_{\ell}$. Let $m, n$ be dual elliptic
pseudoprimes as above and $\ell$ be an Elkies prime. We consider the $\ell$ -
torsion polynomial of $\id{E}_n(C, D)$, which should have $x = \mu \bmod
\id{L}_1$ as an eigenvalue, where $(\ell) = \id{L}_1 \cdot \id{L}_2$ is the
splitting of $\ell$ in $\id{O}(\K)$. If $n$ is prime, there is an Elkies
factor verifying:
\[ F(X^n) - F(g_{x}(X)) \equiv 0 \bmod F(X), \]
where $g_{x}$ is the multiplication polynomial defined in \rf{mulpol} with
respect to $\psi_{\ell}(X)$. For pseudoprime $n$, we let
\begin{eqnarray}
\label{elkf}
h_1(X) & = & X^n \quad \rem \quad \psi_{\ell}(X), \nonumber \\ h_2(X) & = &
\psi_{\ell}(g_x(X)) \quad \rem \quad \psi_{\ell}(X) \quad \hbox{ and } \\ F(X)
& = & \GCD\left(\psi_{\ell}(X), h_1(X) - h_2(X) \right). \nonumber
\end{eqnarray} 
If $x^2 \equiv m \bmod \ell$, then the eigenvalue $x$ is double and we may
discard $\ell$ or use direct factorization, e.g. some variant of the Berlekamp
algorithm \cite{GG}, Chapter V., for finding an Elkies factor.

If $x^2 \not \equiv m \bmod \ell$ and $F(X)$ does not verify the defining
conditions for an Elkies factor, then $n$ must be composite, and the primality
test would stop at this point. Otherwise $F(X)$ is a factor which can be used
in proving existence of an \nth{\ell} elliptic extension.

\section{Applications to Cyclotomy}
We now come to the application of dual elliptic pseudoprimes for the cyclotomy
primality test. A first application of these pseudoprimes was given in
\cite{Mi1} and it took advantage of the Corollary \ref{cdif} and the implied
fourth root order bound \rf{2.4} on the difference between the smallest
eventual divisors of $(m, n)$; this was an improvement on methods for finding
divisors in residue classes, like \cite{Le3}, \cite{CHN}.

By using elliptic extensions and Theorem \ref{thL}, we are in the
more pleasant situation, that trial division may be completely
eliminated in the cyclotomy tests. The particularity of our new
algorithm consists in the inhabitual fact that, for proving primality
of one pseudoprime, it is more efficient to do so for \textit{ two
pseudoprimes } simultaneously. Only this allows, of course, to use the
strong implications of duality.

Suppose that $n$ is a test number like before and a second strong pseudoprime
$m < n$ was found, such that $(m, n)$ are dual elliptic pseudoprimes with
respect to the order $\id{O} \subset \K = \Q(\sqrt{-d})$. We choose some
parameters $s, t$ with $s > 2n^{1/4}$ and $t = \lambda(s)$, the Carmichael
function. Then we find a good $L$ with the algorithm ACE and choose a divisor
$s' | s $ such that for $S = s' \cdot L$, the inequality
\begin{eqnarray}
\label{nsiz}
\left| ( \ m \ \rem \ S \ ) - ( \ n \ \rem \ S \ ) \right| > 2 n^{1/4}
\end{eqnarray}
holds. Next one performs the main stage of the cyclotomy test for $S$, on
\textit{both} $m$ and $n$ and proves the existence of an \nth{L} elliptic
extension by verifying \rf{oddtaue}, \rf{eventaue} in the same working
extensions used for the cyclotomy test. Since $t | \varphi(L)$ and equality is
not necessary, some additional working extensions will in general be
required. Note that in building elliptic Jacobi sums, one has also to check
that the primes involved are not exceptional conductors. If this happens, the
respective $\ell | L$ should be replaced by a new one, keeping the properties
of $L$ valid.

The Theorem 3. implies that there is a prime $p < \sqrt{n}, \ p | n$ and a
dual elliptic $q$ to $p$, which divides $m$. Furthermore, the algorithm ACE
and \rf{fdiv} imply that
\[ | p - q | = \left| ( \ m \ \rem \ S \ ) - ( \ n \ \rem \ S \ ) \right| > 2
n^{1/4}, \] in contradiction with \rf{2.4} and it follows that $m, n$ must be
primes.

We formulate the strategy described above in algorithmic form.

\begin{quotation}
\NI {\bf Algorithm CIDE}{\em ( Cyclotomy Initialized by Dual
Elliptic tests ) }
\vspace*{0.3cm}
\NI Let  $n$ be a strong pseudoprime.
\vspace*{0.1cm}
\begin{itemize}
\item[1.] Find a dual elliptic pseudoprime $m$ to $n$, with respect to an
  order $\id{O} \subset \K = \Q[ \sqrt{-d}]$, by using standard versions of
  ECPP. If none can be found ( in affordable time ), then stop or skip to a
  classical cyclotomy test for $n$.  
\vspace*{0.1cm}
\item[2.] Choose the parameters $s, t$, such that \rf{nsiz} is verified and
 $t = \lambda(s) $ ( \cite{Mi2}). 
\item[3.] Find a good $L$ using algorithm ACE and let $S = L \cdot s'$, where
  $s'$ is the smallest factor of $s$ such that \rf{nsiz} is verified by $S$.
\vspace*{0.1cm}
\item[4.] Construct saturated working extensions of $\ZM{m}, \ZM{n}$ for each
  prime $v | \varphi(L)$. Let $\rg{R}_m, \rg{R}_n$ be their compositum.
\item[5.] Perform in $\rg{R}_m$ respectively $\rg{R}_n$ the Jacobi sum tests
\rf{cpp} necessary for proving the existence of \nth{S} cyclotomic extensions
of $\ZM{n}$ and of $\ZM{m}$.
\item[6.]  Compute the elliptic Jacobi sums related to formulae \rf{oddtaue}
  and \rf{eventaue} for all $\ell | L$ and eventually replace $\ell$ if it is
  an exceptional conductor. 
\item[7.] Perform in $\rg{R}_m$ respectively $\rg{R}_n$ the elliptic Jacobi
  sum tests implied by \rf{oddtaue} and \rf{eventaue}, which are necessary for
  proving the existence of \nth{L} elliptic extensions of $\ZM{n}$ and of
  $\ZM{m}$.
\item[8.] Declare $m$ and $n$ prime if all the above tests are passed
  successfully. 
\end{itemize}
\end{quotation}
\subsection{Run Time}
We split the computations for a CIDE - test for a probable prime $n \in \N$ in
three main stages:
\begin{itemize}
\item[I.]     Find a dual elliptic pseudoprime $m$ to $n$.
\item[II.]    Perform cyclotomy tests for $m, n$.
\item[III.]  Find an $L$ with the ACE algorithm and prove the existence of an
\nth{L} elliptic extension for $m$ and $n$.
\end{itemize}
If the Jacobi sums for the Step II. are computed in essentially linear time,
e.g. by using the algorithm of Ajtai et. al. \cite{AjKuSi}, then Step
II. reduces to the main stage of the cyclotomy test. This stage is polynomial
and takes $O(\log^3(n))$ binary steps \cite{Mi4}. As mentioned above,
heuristic arguments suggest that Step I. also takes cubic time \cite{Mo2},
\cite{GaK}. 

We analyze the run time for the Step III using the heuristics in Fact
\ref{heur}. Let the bound $B$ be defined by \rf{smallp}; the factors of $L$
will be $\ell < B^2$ and their number is $< B$. For each factor, one has to
perform some elliptic Jacobi sum tests, at most $\log(B)$; the degree of the
extensions where the tests are performed is also $ < B^2$. Altogether, using
$B = O(\log \log(n))$, this implies that Step III is performed in
\[ O\left(\log^{2+\varepsilon}(n) \times B^{3 + \varepsilon}\right) = O\left(
\log^{2+\varepsilon}(n)\right) \] binary operations. The Step III. is thus
dominated by steps I and II.  Hence, the run time of the algorithm CIDE is:
\[ O\left(\log(n)^{3 + \varepsilon}\right).\]
 
\begin{remark}
\label{cert}
Using the certification algorithm described in \cite{Mi4}, one can
also provide primality certificates which can be verified in quadratic
time. Note that this time is \textbf{unconditional} and can be
achieved also if no certified Jacobi sum tables are available.
\end{remark}

\section{Conclusions }
Since the summer of $2002$, the theoretical problem of primality
proving is solved: \textit{Primes is in P}, as Agrawal, Kayal and
Saxena laconically put it the title of their magnificent paper
\cite{AKS}.  Apart from the thus closed search for a polynomial time
deterministic algorithm, there is an alternative question concerning
primality proving. Namely: "How large general numbers can
\textbf{currently} be proved on a computer"?

It is a general fact that provable algorithms are different from their
practical versions, which, if they exist, may lose some or many of the
theoretical advantages, but work conveniently in practice. Thus, the algorithm
of Goldwasser and Kilian \cite{GoKi,GoKi99} has been proved to terminate in
random polynomial time for all but an exponentially thin set of inputs; it has
hardly ever been implemented, for complexity reasons mentioned in the
introduction. In exchange, the ideas of Atkin \cite{AtMo} led to the current
wide spread version of ECPP \cite{ECPP}, which works very well in practice.
As already mentioned, the choice of the fields of complex multiplication is in
this version such that no \textit{proof} of polynomial time termination is
known; however, the algorithm works very stably in practice and heuristic
argument brought in \cite{GaK} explain this fact.

The situation is even more bizarre with the cyclotomy test: from the
complexity theoretical point of view, it should even not be taken into
consideration, since it is over-polynomial. For the range of primes which are
currently affordable for computer proofs, it works very efficiently. A
fortiori, the combination of cyclotomy and elliptic curves provided by CIDE
has good reasons to be the medium term provider of largest primality proofs
and the generation of certificates which can be verified in quadratic time, as
observed in Remark \ref{cert}, is also an appealing novelty. Furthermore, the
algorithm has random cubic run-time, based on the heuristics of \cite{GaK} and
the ones in Fact \ref{heur}.

Finally, the test of Agrawal, Kayal and Saxena has, for computer
implementation, a serious space problem. Even the nice idea of Berrizbeitia
\cite{Be}, \cite{Bern, AM} which brings an important run - time
improvement\footnote{see various forms of generalizations in \cite{AM},
\cite{Bern}, \cite{Cheng03}, \cite{Mo1},}, does not remove this problem. It is
not likely that primes larger then $500$ decimal digits, say, will be proved
in the near future with any variation of the AKS algorithm, \textit{unless new
ideas are found, for solving the space problem}.

In conclusion, it is a mathematically appealing and relevant goal, to seek for
an efficient variant of AKS, while on the side of CPP, the construction of
Jacobi sums remains a small problem, which is interesting per se. The
algorithm of Ajtai, Kumar and Sivakumar yields however a random polynomial
solution which is satisfactorily in theory, while the LLL and PARI approaches
may solve the practical problem for conceivable applications during the next
years or even decades.  

\bibliographystyle{abbrv} 
\bibliography{cide2}

\begin{thebibliography}{10}

\bibitem{APR}
L.~Adleman and R.~R. C.~Pomerance.
\newblock On distinguishing prime numbers from composite numbers.
\newblock {\em Ann. Math.}, 117:173--206, 1983.

\bibitem{AKS}
M.~Agrawal, N.~Kayal, and N.~Saxena.
\newblock Primes is in p.
\newblock {\em Annals of Math.}, pages 781--793, 2004.

\bibitem{AjKuSi}
M.~Ajtai, R.~Kumar, and D.~Sivakumar.
\newblock A sieve algorithm for the shortest vector problem.
\newblock In A.~Pres, editor, {\em Proceedings of the 33-rd Symposium on Theory
  of Computing (STOC)}, pages 601--610, 2001.

\bibitem{AtMo}
A.~Atkin and F.Morain.
\newblock Elliptic curves and primality proving.
\newblock {\em Math. Comp.}, 61:29--68, 1993.

\bibitem{AM}
R.~Avanzi and P.~Mih\u{a}ilescu.
\newblock Efficient ``quasi''- deterministic primality test improving aks.
\newblock submitted.

\bibitem{Bern}
D.~J. Bernstein.
\newblock Proving primality in essentially quartic random time.
\newblock {\em Math.\ Comp.}, 76(257):389--403, January 2007.

\bibitem{Be}
P.~Berrizbeitia.
\newblock Sharpening 'primes in p' for a large family of numbers.
\newblock {\em Math.\ Comp.}, 74:2043--59, 2005.

\bibitem{BH}
W.~Bosma and M.~der Hulst.
\newblock {\em Primality proving with cyclotomy}.
\newblock PhD thesis, Universiteit van Amsterdam, 1990.

\bibitem{BMSS}
A.~Bostan, F.~Morain, B.~Salvy, and R.~Schost.
\newblock Fast algorithms for computing isogenies between elliptic curves.
\newblock {\em Math. Comp.}, 2007.
\newblock to appear.

\bibitem{Cheng03}
Q.~Cheng.
\newblock Primality proving via one round in {ECPP} and one iteration in {AKS}.
\newblock In D.~Boneh, editor, {\em Advances in Cryptology -- CRYPTO 2003},
  number 2729 in Lecture Notes in Comput. Sci., pages 338--348. Springer
  Verlag, 2003.

\bibitem{CoLe}
H.~Cohen and H.W.{Lenstra Jr.}
\newblock Primality testing and jacobi sums.
\newblock {\em Math. Comp.}, 48:297--330, 1984.

\bibitem{CHN}
D.~Coppersmith, N.~Howgrave-Graham, and S.~V. Nagaraj.
\newblock Divisors in residue classes, constructively.

\bibitem{Cox}
D.~A. Cox.
\newblock {\em Primes of the Form $x^2 + n y^2$}.
\newblock Wiley \& Sons, 1989.

\bibitem{CP}
R.~Crandall and C.~Pomerance.
\newblock {\em Prime Numbers - A Computational Perspective}.
\newblock Springer, 2002.

\bibitem{EnMo02}
A.~Enge and F.~Morain.
\newblock Comparing invariants for class fields of imaginary quadratic fields.
\newblock In C.~Fieker and D.~R. Kohel, editors, {\em Algorithmic Number
  Theory}, volume 2369 of {\em Lecture Notes in Comput. Sci.}, pages 252--266.
  Springer-Verlag, 2002.
\newblock Proceedings of the 5th International Symposium, ANTS-V, Sydney,
  Australia, July 2002.

\bibitem{EnMo03}
A.~Enge and F.~Morain.
\newblock Fast decomposition of polynomials with known {G}alois group.
\newblock In M.~Fossorier, T.~H{\o}holdt, and A.~Poli, editors, {\em Applied
  Algebra, Algebraic Algorithms and Error-Correcting Codes}, volume 2643 of
  {\em Lecture Notes in Comput. Sci.}, pages 254--264, 2003.
\newblock Proceedings, 15th International Symposium, AAECC-15, Toulouse,
  France, May 2003,.

\bibitem{GaK}
S.~D. Galbraith and J.~McKee.
\newblock The probabilitythat the number of points on an elliptic curve over a
  finite field is prime.
\newblock {\em Journal of London Mathematical Society}, 62:671--684, 2000.

\bibitem{GoKi}
S.~Goldwasser and J.~Kilian.
\newblock Almost all primes can be quickly certified.
\newblock In {\em Proc. 18-th Annual ACM Symp. on Theory of Computing}, pages
  316--329, 1986.

\bibitem{GoKi99}
S.~Goldwasser and J.~Kilian.
\newblock Primality testing using elliptic curves.
\newblock {\em Journal of the ACM}, 46(4):450--472, July 1999.

\bibitem{IR}
K.~Ireland and M.~Rosen.
\newblock {\em A Classical Introduction to Modern Number Theory}, volume~84 of
  {\em Springer Graduate Texts in Mathematics}.
\newblock Springe, 1990.
\newblock Second Edition.

\bibitem{Le1}
H.~W. {Lenstra Jr.}
\newblock Primality testing algorithms (after adleman, rumely and williams).
\newblock In {\em Seminaire Bourbaki \# 576}, volume 901 of {\em Lectures Notes
  in Mathematics}, pages 243--258, 1981.

\bibitem{Le3}
H.~W. {Lenstra Jr}.
\newblock Divisors in residue classes.
\newblock {\em Math. Comp.}, pages 331--334, 1984.

\bibitem{Le2}
H.~W. {Lenstra Jr.}
\newblock {\em Galois Theory and Primality Testing}, chapter~12, pages
  169--189.
\newblock Number 1142 in Lecture Notes in Mathematics. Springer Verlag, 1985.

\bibitem{Le4}
H.~W. {Lenstra Jr}, October 2002.
\newblock Seminar lectures in Leiden and Oberwolfach.

\bibitem{LP}
H.~W. {Lenstra, Jr.} and C.~Pomerance.
\newblock Primality testing with gaussian periods.

\bibitem{Mi1}
P.~Mih\u{a}ilescu.
\newblock {\em Cyclotomy of Rings \& Primality Testing}.
\newblock PhD thesis, ETH Z\"urich, 1997.

\bibitem{Mi2}
P.~Mih\u{a}ilescu.
\newblock Cyclotomy primality proving - recent developments.
\newblock In {\em Proceedings of the Third International Symposium ANTS III,
  Portland, Oregon}, volume 1423 of {\em Lecture Notes in Computer Science},
  pages 95--111, 1998.

\bibitem{Mi4}
P.~Mih\u{a}ilescu.
\newblock Cyclotomy primality proofs and their certificates.
\newblock {\em Mathematica Gottingensis}, 2006.

\bibitem{MMS}
P.~Mih\u{a}ilescu, F.~Morain, and E.~Schost.
\newblock Computing the eigenvalue in the schoof - elkies - atkin algorithm
  using abelian lifts.
\newblock In {\em ISSAC}, 2007.

\bibitem{MV}
P.~Mih\u{a}ilescu and V.~Vuletescu.
\newblock Elliptic gauss sums and counting points on curves, 2007.

\bibitem{Mo4}
F.~Morain.
\newblock private communication.

\bibitem{ECPP}
F.~Morain.
\newblock Site for downloading the elliptic curve primality test software of f
  .morain.".

\bibitem{Mo3}
F.~Morain.
\newblock Calcul du nombre de points sur une courbe elliptique dans un corps
  fini~: aspects algorithmiques.
\newblock {\em Journal de Th\'eorie des Nombres, Bordeaux}, 7:255--282, 1995.

\bibitem{Mo}
F.~Morain.
\newblock Primality proving using elliptic curves: An update.
\newblock In {\em Proceedings of the Third International Symposium ANTS III,
  Portland, Oregon}, volume 1423 of {\em Lecture Notes in Computer Science},
  pages 111--127, 1998.

\bibitem{Mo1}
F.~Morain.
\newblock La primalit\'e en temps polynomial [d'apr\`es adleman, huang;
  agrawal, kayal, saxena].
\newblock In {\em Seminaire Bourbaki 55-\`eme ann\'ee}, number 917 in Lectures
  Notes in Mathematics, 2002-2003.

\bibitem{Mo2}
F.~Morain.
\newblock Implementing the asymptotically fast version of elliptic curve
  primality proving algorithm.
\newblock {\em Math. Comp}, 76(257), January 2007.

\bibitem{Pi}
R.~Pinch.
\newblock Galois module structure of elliptic functions.
\newblock In O.~U. Press, editor, {\em Computers in mathematical research
  (Cardiff, 1986)}, number~14 in Inst. Math. Appl. Conf., pages 69--91, 1988.
\newblock New Series.

\bibitem{Sch}
R.~Schoof.
\newblock Counting points on elliptic points over finite fields.
\newblock {\em J. Th. Nombr. Bordeaux}, 7:363--397, 1995.

\bibitem{Si}
J.~Silverman.
\newblock {\em The Arithmetic of Elliptic Curves}, volume 106 of {\em Graduate
  Texts in Mathematics}.
\newblock Springer, 1996.

\bibitem{GG}
J.~{von zur Gathen} and J.~Gerhardt.
\newblock {\em Modern Computer Algebra}.
\newblock Cambridge University Press, 2-nd ed. edition, 2000.

\bibitem{W}
L.~Washington.
\newblock {\em Elliptic Curves - Number Theory and Cryptography}.
\newblock Chapman and Hall/CRC, 2003.

\end{thebibliography}

\end{document}